\documentclass[aip, prb, numerical,preprint,nofootinbib]{revtex4-1}
\usepackage[english,activeacute]{babel}
\usepackage{amssymb,amsmath,alltt}
\usepackage[utf8]{inputenc}
\usepackage{amsfonts}
\usepackage{graphicx}
\usepackage{multirow}

\usepackage{amsthm}
\newtheorem{definition}{Definition}
\newtheorem*{remark}{Remark}
\newtheorem{theorem}{Theorem}
\newtheorem{algorithm}{Algorithm}
\newtheorem{proposition}[theorem]{Proposition}

\def\torus{{\mathbb T}}

\newcommand{\bq}{\begin{equation}}
\newcommand{\eq}{\end{equation}}

\newcommand{\be}{\begin{eqnarray}}
\newcommand{\ee}{\end{eqnarray}}

\begin{document}
\title{Global transport in a nonautonomous standard map}
\author{R. Calleja}
\email{calleja@mym.iimas.unam.mx}
\affiliation{IIMAS-UNAM, Mexico D.F., 04510, Mexico}
\author{D. del-Castillo-Negrete}
\email{delcastillo@ornl.gov}
\affiliation{Oak Ridge National Laboratory, Oak Ridge, Tennessee, 37831-8071, USA}
\author{D. Martínez-del-Río}
\email{dmr@mym.iimas.unam.mx}
\affiliation{IIMAS-UNAM, Mexico D.F., 04510, Mexico}
\author{A. Olvera}
\email{aoc@mym.iimas.unam.mx}
\affiliation{IIMAS-UNAM, Mexico D.F., 04510, Mexico}

\begin{abstract}
A non-autonomous version of the standard map 
with a periodic variation of the parameter
is  introduced and studied.  
Symmetry properties in the variables and parameters of the map are found and used to find relations between rotation numbers of invariant sets of the autonomous realization of the period-two case of the map. 
The role of the nonautonomous dynamics on period-one orbits, stability and bifurcation is studied.
 The critical boundaries for the global transport and the destruction of  invariant circles with fixed rotation number are studied in detail using direct computation and a continuation method. In the case of global transport, the critical boundary  has a particular symmetrical horn shape. The results are contrasted with similar calculations found in the literature.
\end{abstract}

\maketitle

{\bf
A nonautonomous version of the standard map is presented, using an autonomous realization of the map, the existence of invariant tori (and associated global transport problem) as a function of the parameters is studied analytically, numerically, and with {a continuation} method. 
The study of the map symmetries allows to reduce the parameter space and infer the existence of invariant circles.
The results show a particular symmetrical horn shaped critical boundary for global transport in the parameter space that contains all the critical boundaries for invariant tori with fixed rotation numbers calculated by the parametrization method. 
}
%|||||||||||||||||||||||||||||||||||||||||||||||||||||
\section{Introduction}
%|||||||||||||||||||||||||||||||||||||||||||||||||||||
Two dimensional twist maps have been extensively studied within discrete dynamical systems \cite{Meiss92}. However, studies of non-autonomous maps seem to be scarce, even for two dimensional cases. Here the term \emph{non-autonomous} is used to describe any map $T:\mathtt{x}_n \mapsto \mathtt{x}_{n+1}$  that have explicit dependence on the iteration number: $\mathtt{x}_{n+1}=f(\mathtt{x}_n,n)$. The present work illustrates that a non-autonomous periodic variation in a parameter can give rise to unexpected behavior in the well known standard map.

Our motivation to define a \emph{non-autonomous standard map} (NASM) comes from a
simplified model of self-consistent transport in 
marginally stable systems including vorticity mixing in strong shear flows and electron dynamics in plasmas, known as the \emph{self-consistent map model} or \emph{single wave map model}.
The model, originally presented in Ref.~\onlinecite{del_castillo_2000} and later on studied in 
Ref.~\onlinecite{Bofetta03, Carbajal12, Coreo15}
consists of a set of $N$ mean-field coupled standard-like area-preserving twist maps in which the amplitude and phase of the perturbation (rather than being constant) are dynamical variables.
The numerical simulations of this map for a particular kind of initial conditions shows the existence of coherent structures, which lead in Ref.~\onlinecite{Coreo15} to introduce a two dimensional map model that mimics the asymptotic behavior of the original map and sheds light on the effects of an external oscillatory field over an uncoupled system, the NASM. 

This kind of maps has been studied before in the literature \cite{canadel15,tompaidis96,simo16} following different approaches. Ref. \onlinecite{canadel15} considered families of non-autonomous maps that  converge to a known autonomous map, with known invariant tori, and  
 studied the convergence to \emph{asymptotic invariant tori}. 
 The \emph{rotating standard map} in Ref. \onlinecite{tompaidis96} and the \emph{driven standard map} in Ref. \onlinecite{simo16} can be considered more general cases of the NASM defined here, although with different aims: 
 a search of the bifurcation space where the two dimensional tori cease to exist and regions of stability of invariant curves in the parameter space, respectively.

The goal of the present paper is to give a detailed characterization of the \emph{non-autonomous standard map} (NASM) defined in Ref.~\onlinecite{Coreo15}  and to determine  
conditions for the existence of global transport. 
To do this, we start by defining the map in Sec. \ref{sec:MapDefinition} along with an autonomous realization of the map and some general properties. 
 In Sec. \ref{sec:Symmetries},
  we study the map's symmetries in the coordinates and parameters spaces, and find relations between rotation numbers of invariant circles for symmetric values of the parameters. 
In section \ref{sec:Periodic}, we study 
the periodic orbits 
of the NASM. 
In Sec. \ref{sec:KnownLimits}, we study the known reduced cases of the map and give a criteria to search numerically for the threshold of global transport. 
In Sec. \ref{sec:Continuation}, we  search again for transport barriers as function of the parameters, using the \emph{parametrization method}, {applied} on particular interesting rotation numbers.  
 In Sec. \ref{sec:Discussion}, we compare these results with  previous works in Refs. \onlinecite{canadel15, tompaidis96,simo16} and propose some conjectures
and
present the conclusions.

%|||||||||||||||||||||||||||||||||||||||||||||||||||||
\section{Map definition.}
%|||||||||||||||||||||||||||||||||||||||||||||||||||||
\label{sec:MapDefinition}
Our starting point is the nonautonomous standard map (NASM), defined as
\begin{subequations}
\label{nas_map}
\begin{align}
   \hat{x}_{n+1} &= \hat{x}_n + \hat{y}_{n+1} 
   \;\;\;\;\;\;\;\;\;\;\; \mathrm{mod} \;\, (1), \label{nas_map1}\\
   \hat{y}_{n+1}&= \hat{y}_n + \frac{\kappa_n}{2\pi} \sin(2\pi \hat{x}_n)\,. \label{nas_map2} 
\end{align}
\end{subequations}
where $\kappa_n$ is a function of $n$.

Motivated by the asymptotic dynamics of self-consistent coupled maps\cite{del_castillo_2000,Coreo15}, we focus on the case when $\kappa_n$ is a periodic function. As a first step, we  consider  a ``triangular wave" periodic dependence in which 
$\kappa_n$ can only take two values:
\begin{equation}
  \kappa_n = \left\lbrace \begin{array}{cl}
    \kappa_1 & \;\;\;\; {\rm if} \; n \; {\rm is}\;{\rm odd}, \\ & \\
    \kappa_2 & \;\;\;\; {\rm if} \; n \; {\rm is}\;{\rm even}. \\
\end{array} \right.
\label{kap}
\end{equation}

\begin{figure}[h!]
     \centering
     \includegraphics[width=\linewidth, trim={1 1 0 0}]{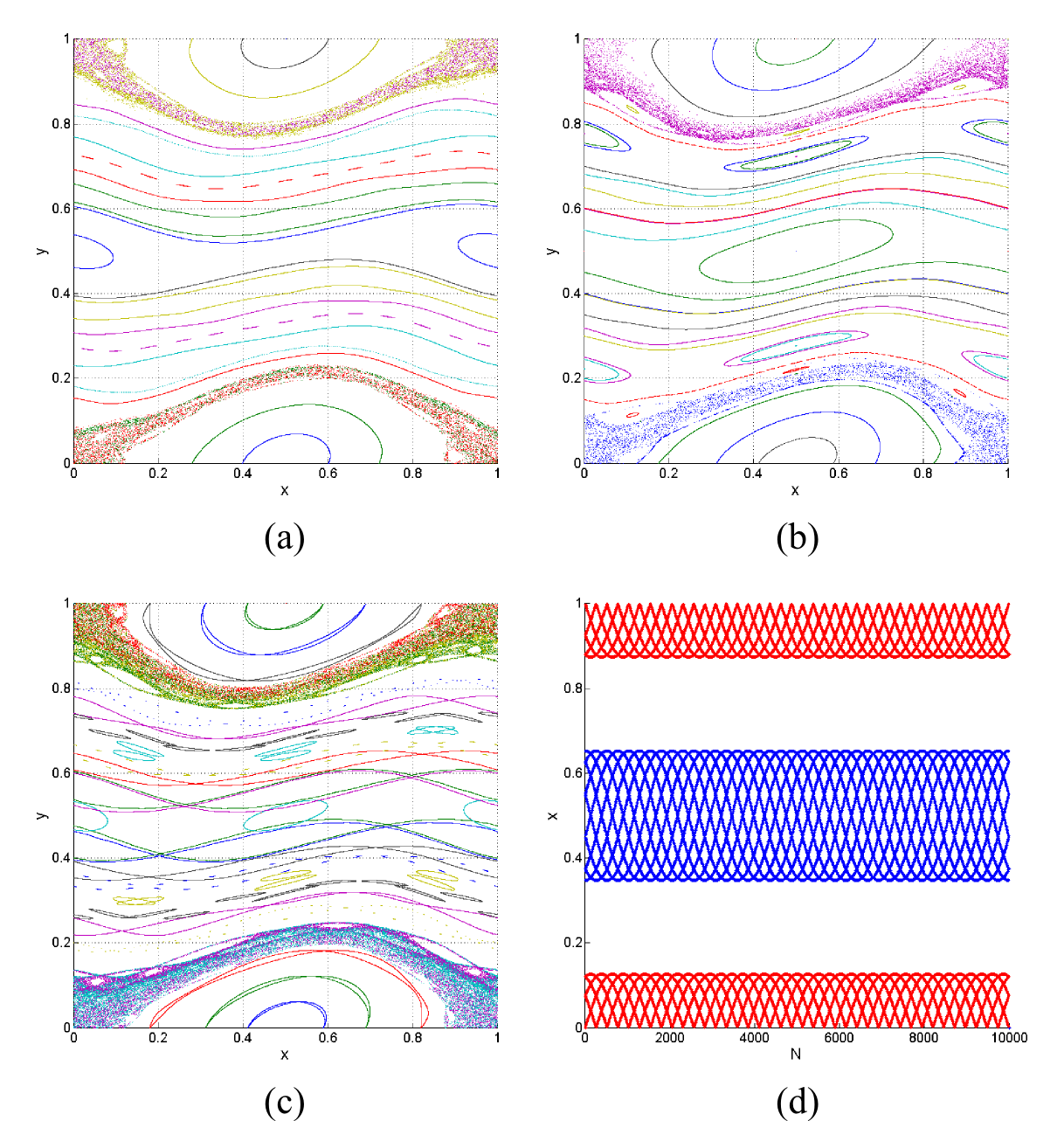}
     \caption{Comparison between the phase space of maps (\ref{mapaut}) and (\ref{nas_map}). (a) and (b) correspond to (\ref{mapaut})  with $(\kappa_1,\kappa_2)$: $(0.35,0.5)$  and  $(0.5,0.35)$ respectively, and (c) corresponds to map (\ref{nas_map}) with the same parameters. 
     (d)  shows to the time series of the initial condition $z_0=(0,0.55)$, where blue (red) corresponds to the even (odd) iterates of map (\ref{nas_map}).}
     \label{maps_phase_space}
 \end{figure}

 We define a new map, $\mathcal{T}_{\kappa_1 \kappa_2}$, such that its iterates $(x_n,y_n)$ $n=1,2,...$ coincide with the even iterations of (\ref{nas_map}),  i.e.,
\begin{equation}
(x_n,y_n) = (\hat{x}_{2n},\hat{x}_{2n})\,.
\label{hat_no_hat}
\end{equation}
 By construction, the map $\mathcal{T}_{\kappa_1 \kappa_2}$ is autonomous  
and can be written as,
\begin{subequations}
\label{mapaut}
\begin{align}
   {x}_{n+1} &= {x}_n + 2{y}_n +
      \mathcal{F}_1({x}_n,{y}_n;\kappa_1,\kappa_2)
       \;\;\;\;\;\; \mathrm{mod} \;\, (1) \label{mapaut1}\\
   {y}_{n+1}&= {y}_n +
             \mathcal{F}_2({x}_n,{y}_n;\kappa_1,\kappa_2) \label{mapaut2} 
\end{align}
\end{subequations}
 where  
the functions $\mathcal{F}_1$ and $\mathcal{F}_2$ are defined as,
\begin{subequations}
\label{EfeGe}
\begin{align}
     &\mathcal{F}_1({x},{y};\kappa_1,\kappa_2) =  \frac{\kappa_1}{2\pi} \sin(2\pi {x})+ \mathcal{F}_2
      \, , \label{Efe}\\
      &\mathcal{F}_2({x},{y};\kappa_1,\kappa_2)= \frac{\kappa_1}{2\pi} \sin(2\pi {x}) \nonumber \\ 
       & \phantom{OOOO} + \frac{\kappa_2}{2\pi} \sin \left\lbrace 2\pi \left[ {x} + {y} + \frac{\kappa_1}{2\pi} \sin (2 \pi {x}) \right]\right\rbrace \, . \label{Ge}
\end{align}
\end{subequations}
 As Fig.~\ref{maps_phase_space}(c) shows, due to its non autonomous nature, the map in (\ref{nas_map}) exhibit intersection of orbits, something that as Fig.~\ref{maps_phase_space}(a) and Fig.~\ref{maps_phase_space}(b) illustrates, never happens in autonomous maps. Fig.~\ref{maps_phase_space}(d) shows, the non autonomous dynamics of an initial condition inside an \emph{island} of (\ref{nas_map}) alternates between $(0,1/2)$ elliptic point of $\mathcal{T}_{\kappa_1 \kappa_2}$ and $(1/2,1/2)$ elliptic point of $\mathcal{T}_{\kappa_2 \kappa_1}$.

%////////////////////////////
%////////////////////////////

The map in Eq.~(\ref{mapaut}) is equivalent to the composition of two standard maps, i.e., 
$\mathcal{T}_{\kappa_1 \kappa_2}\equiv\mathcal{S}_{\kappa_2}\circ\mathcal{S}_{\kappa_1}$, where 
$\mathcal{S}_{\epsilon}$ denotes the  standard map
with perturbation parameter $\epsilon$,
\begin{equation}
 \mathcal{S}_\epsilon 
 \left( \begin{array}{c}  x\\ y \end{array} \right)  =
    \left(\begin{array}{c}
   {x}+ y+ \frac{\epsilon}{2\pi} \sin(2 \pi x) \\ 
 y+ \frac{\epsilon}{2\pi} \sin(2 \pi x) 
\end{array} \right)\,.
\label{stdmap}
\end{equation}
with $(x,y)\in\mathbb{T}\times\mathbb{R}$.
It is straightforward to show that $\mathcal{T}_{\kappa_1 \kappa_2}$ and $\mathcal{T}_{\kappa_2 \kappa_1}$ are  diffeomorphic, since,
\bq
\label{sim_M}
\mathcal{S}_{\kappa_2}^{-1} \circ {\cal T}_{\kappa_1 \kappa_2} ^n \, \circ \mathcal{S}_{\kappa_2} =   {\cal T}_{\kappa_2 \kappa_1}  ^n \,.
\eq
where $\mathcal{S}^{-1}_{\epsilon}$ denotes the inverse of the standard map,
\begin{equation}
 \mathcal{S}_\epsilon^{-1} 
 \left( \begin{array}{c}  x\\ y \end{array} \right)  =
    \left(\begin{array}{c}
   {x} - y \\
 y - \frac{\epsilon}{2\pi} \sin(2 \pi (x-y)) 
\end{array} \right)\,.
\label{stdmap_inv}
\end{equation}

The \emph{twist} of $\mathcal{T}_{\kappa_1 \kappa_2}$ is given by,
\begin{eqnarray}
   \frac{\partial x_{n+1}}{\partial y_n}\vert_{x_n}&=&
    2+\kappa_2\,\cos\left(2\pi \left(x_n+y_n+\frac{\kappa_1}{2\pi}\,\sin(2\pi x_n) \right)\right).
   \label{twist_cond} \nonumber \\
   &&
\end{eqnarray}
and the twist condition $\frac{\partial x_{n+1}}{\partial y_n}|_{x_n}>0$, is satisfied in the whole cylinder when  $|\kappa_2|<2$. 
Based on this, we
define the \emph{twist region} in the parameter space as the square: $\{ (\kappa_1,\kappa_2) : |\kappa_i| < 2, i=1,2\}$.
It must be noted that, if the map is twist for small values of the parameters, 
the 
Moser's Twist Theorem \cite{moser62} guarantees the existence of invariant circles   well ordered  with the rotation number.

%===================================================
Let $\Omega = dy\wedge dx$ be a symplectic form on the cylinder, then it is clear
that $\mathcal{T}_{\kappa_1 \kappa_2}$ is symplectic with respect to $\Omega$ since each one of 
the standard maps is symplectic. In particular, 
\begin{equation}
(\mathcal{T}_{\kappa_1 \kappa_2})^* \Omega = (\mathcal{S}_{\kappa_2}\circ\mathcal{S}_{\kappa_1})^* \Omega
 = \mathcal{S}_{\kappa_1}^*\mathcal{S}_{\kappa_2}^* \Omega = \Omega\,,
\label{symplectic}
\end{equation}
{where $(\mathsf{F})^*$ is the pullback via a function $\mathsf{F}$}.
We also note that since both
$\mathcal{S}_{\kappa_2}$ and $\mathcal{S}_{\kappa_2}$ are exact, then $\mathcal{T}_{\kappa_1 \kappa_2}$ is also 
exact. Notice that if $ d \alpha = \Omega$, then
\[
\mathcal{S}_{\kappa_i}^* \alpha - \alpha = d P_i \qquad i=1, 2
\]
where $P_i$ is the generating functions of $\mathcal{S}_{\kappa_i}$.
The simple computation is as follows (see Ref.~\onlinecite{Greene-14}). 
\begin{eqnarray}\label{exact}
\mathcal{T}_{\kappa_1 \kappa_2}^* \alpha &=& (\mathcal{S}_{\kappa_2}\circ\mathcal{S}_{\kappa_1})^*\alpha =
\mathcal{S}_{\kappa_1}^* \mathcal{S}_{\kappa_1}^*\alpha = \mathcal{S}_{\kappa_1}^* (\alpha + d  P_2) \nonumber \\
       &=& \alpha + d( P_1 +  \mathcal{S}_{\kappa_1}^* P_2)
\end{eqnarray}
Therefore $\mathcal{T}_{\kappa_1 \kappa_2}$ is exact with generating function $P_1 +  \mathcal{S}_{\kappa_1}^* P_2$.
%=====================================================

%||||||||||||||||||||||||||||||||||||||||||||||||||
\section{Symmetries.}
%||||||||||||||||||||||||||||||||||||||||||||||||||
\label{sec:Symmetries}
 The symmetries of the NASM, 
help to reduce different possible 
cases and  to infer 
the existence of invariant circles and their rotation numbers 
 from the existence of invariant circle with different parameter values.
 Since the autonomous realization of the NASM in (\ref{mapaut}) is the composition of two standard maps (\ref{stdmap}),  we expect to have 
 the symmetries of the standard map. 
However there are other  
 symmetries unique to $\mathcal{T}_{\kappa_1 \kappa_2}$. 

The functions $\mathcal{F}_i$ in Eqs.~(\ref{Efe}) and (\ref{Ge}) have the following symmetries.
\begin{eqnarray}
   P1.\,\,\,  &&
       \mathcal{F}_i(-x,-y;\kappa_1,\kappa_2)=-\mathcal{F}_i(x,y;\kappa_1,\kappa_2).
       \label{sym_1} \\
   P2.\,\,\, &&
      \mathcal{F}_i(x+n,y+m;\kappa_1,\kappa_2)=\mathcal{F}_i(x,y;\kappa_1,\kappa_2), \,\,\, \nonumber \\ 
           && \mathrm{for}\,\, n,m\in \mathbb{Z}. \label{sym_2} \\
   P3.\,\,\, &&
       \mathcal{F}_i(x,y; \kappa_1;\kappa_2)=\mathcal{F}_i(x+1/2,y; -\kappa_1,-\kappa_2).   
      \label{sym_3}   \\
   P4.\,\,\, &&
      \mathcal{F}_i(x,y; \kappa_1;\kappa_2)=\mathcal{F}_i(x,y+1/2; \kappa_1,-\kappa_2).
       \label{sym_4}  
\end{eqnarray}
Whereas Eqs.~(\ref{sym_1}) and (\ref{sym_2}) are coordinate symmetries directly inherent{ed} from the standard  map,
Eqs.~(\ref{sym_3}) and (\ref{sym_4}) are unique symmetries of the non-autonomous map that involve both coordinates and  parameter transformations.
From,  (\ref{sym_3}) and (\ref{sym_4}) it follows that,
 \begin{eqnarray}
   &&
       \mathcal{F}_i(x,y; \kappa_1;\kappa_2) =  \mathcal{F}_i(x+1/2,y+1/2; -\kappa_1; \kappa_2). \phantom{O}
      \label{sym_5}    
\end{eqnarray}

 In the  remaining of this section, 
 we denote by $(x,y)$ 
 the variables of the \emph{lift} of map (\ref{mapaut}). In other words,  $x\in\mathbb{R}$.

%--------------------------------------------------
\subsection{Orbit symmetries.}
\begin{figure}[h!]
     \centering
     \includegraphics[width=\linewidth]{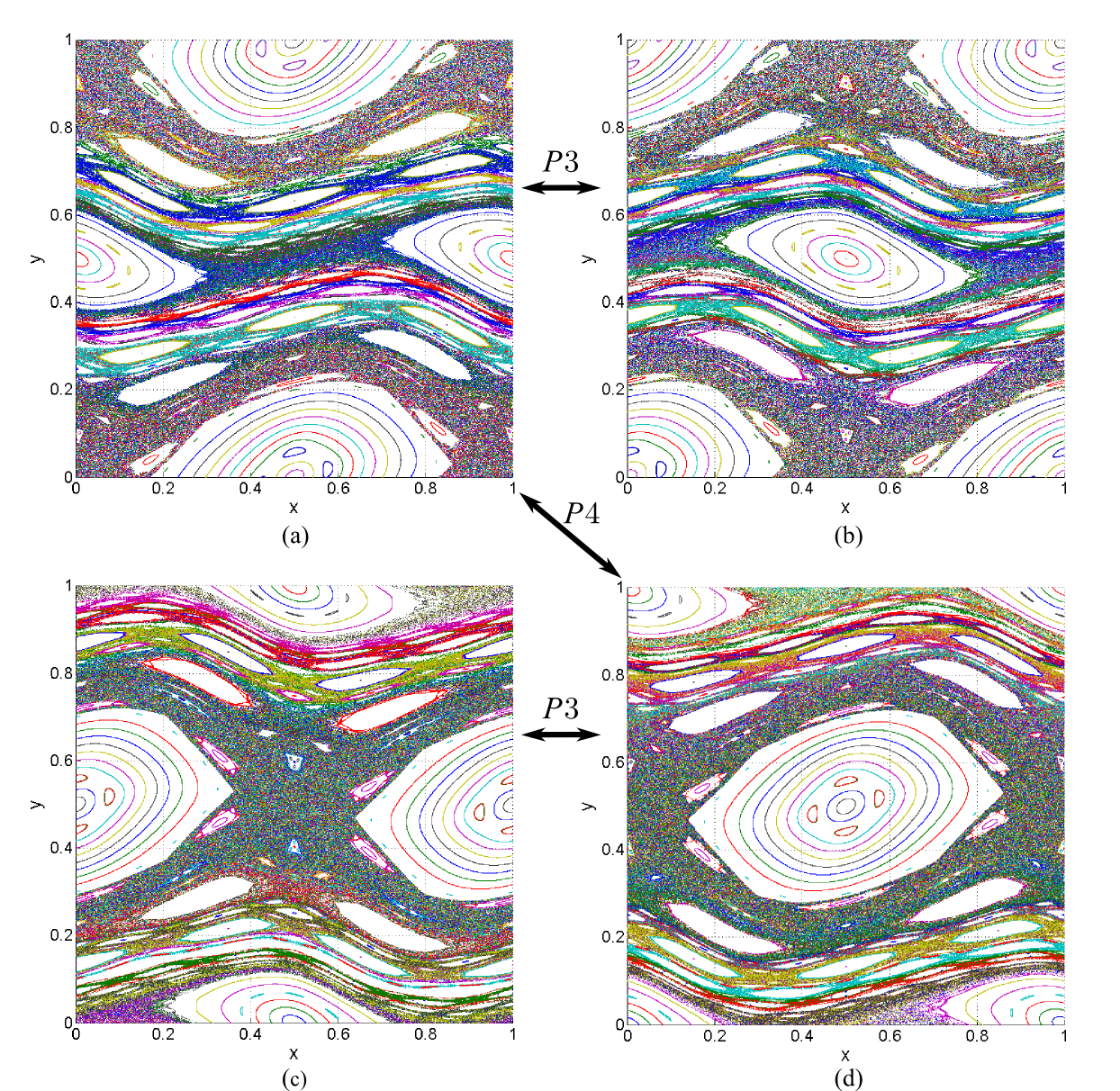} 
     \caption{Phase space of the map (\ref{mapaut}) for different values of the parameters $(\kappa_1,\kappa_2)$: (a) $(0.5,0.7)$ ,(b) $(-0.5,-0.7)$, (c) $(-0.5,0.7)$ and (d) $(0.5,-0.7)$, related by the symmetries $P3$ and $P4$, in Eq.~(\ref{sym_3}) and (\ref{sym_4}) respectively.}
     \label{cell_sym}
 \end{figure}

%--------------------------------------------------
Let $x_n(x_0,y_0; \kappa_1,\kappa_2)$ and $y_n(x_0,y_0; \kappa_1,\kappa_2)$ denote the $x$ and $y$ coordinates of the $n$-th iterate of the  NASM with parameters $(\kappa_1,\kappa_2)$ and initial condition $(x_0,y_0)$.
Then, from the properties of $\mathcal{F}_i$, it follows that the orbits exhibit  the following symmetries (see Fig.~\ref{cell_sym}). 

% ========================================
\begin{enumerate}
  %!!!!!!!!!!!!!!!!! - 1 - !!!!!!!!!!!!!!!!!!!!
  \item Coordinate reflection [from Eq.~(\ref{sym_1})],
       \begin{subequations}
       \label{refl}
       \begin{align}
     {x}_n(-{x}_0,-{y}_0; \kappa_1,\kappa_2) &= - {x}_n({x}_0,{y}_0;\kappa_1,\kappa_2) \, ,
            \label{refl_p1}\\
   {y}_n(-{x}_0,-{y}_0; \kappa_1,\kappa_2) &= - {y}_n({x}_0,{y}_0;\kappa_1,\kappa_2)  \,
           \label{refl_p2} . 
       \end{align}
       \end{subequations}
       
  %!!!!!!!!!!!!!!!!! - 2 - !!!!!!!!!!!!!!!!!!!!  
  \item Coordinate translation [from Eq.~(\ref{sym_2})],
     \begin{subequations}
     \label{trans}
     \begin{align}
      {x}_n({x}_0+r,{y}_0+s; \kappa_1,\kappa_2) &=
         {x}_n({x}_0,{y}_0;\kappa_1,\kappa_2) +  r +2ns  \,,
        \label{trans_p1}\\ 
      {y}_n({x}_0+r,{y}_0+s; \kappa_1,\kappa_2) &= 
         {y}_n({x}_0,{y}_0;\kappa_1,\kappa_2) +s \, ,
          \label{trans_p2}
       \end{align}
       \end{subequations}
   for $r,s\in \mathbb{Z}$. 
 
 %!!!!!!!!!!!!!!!!!! - 3 - !!!!!!!!!!!!!!!!!!!!!!!
  \item Coordinate translation and reflection [from Eqs.~(\ref{refl})-(\ref{trans})], 
     \begin{subequations}
     \label{int}
     \begin{align}  
      {x}_n(1-{x}_0,1-{y}_0; \kappa_1,\kappa_2) &= 
        1- {x}_n({x}_0,{y}_0;\kappa_1,\kappa_2) + 2n \, ,
        \label{int_p1}\\ 
      {y}_n(1-{x}_0,1-{y}_0; \kappa_1,\kappa_2) &= 
         1- {y}_n({x}_0,{y}_0;\kappa_1,\kappa_2)  \, 
          \label{int_p2}. 
       \end{align}
       \end{subequations}

  %!!!!!!!!!!!!!!!!!! - 4 - !!!!!!!!!!!!!!!!!!!!!!!
   \item Coordinate translation  and parameter reflexion [from Eq.~(\ref{sym_3})],
     \begin{subequations}
     \label{shift3}
     \begin{align}
      {x}_n({x}_0+1/2,{y}_0; -\kappa_1,-\kappa_2) &= 
          {x}_n({x}_0,{y}_0;\kappa_1,\kappa_2) +1/2 \,, 
      \label{shift31}\\
      {y}_n({x}_0+1/2,{y}_0; -\kappa_1,-\kappa_2) &= 
        {y}_n({x}_0,{y}_0;\kappa_1,\kappa_2)  \, 
      \label{shift32} . 
     \end{align}
     \end{subequations}

  %!!!!!!!!!!!!!!!!!! - 5 - !!!!!!!!!!!!!!!!!!!!!!!
  \item Coordinate translation  and parameter reflexion [from Eq.~\ref{sym_4}],
     \begin{subequations}
     \label{shift2}
     \begin{align}
       {x}_n({x}_0,{y}_0+1/2; \kappa_1,-\kappa_2) &= 
           {x}_n({x}_0,{y}_0;\kappa_1,\kappa_2) + n \, ,    
       \label{shift21}\\
      {y}_n({x}_0,{y}_0+1/2; \kappa_1,-\kappa_2) &= 
           {y}_n({x}_0,{y}_0;\kappa_1,\kappa_2) +1/2 \, 
      \label{shift22}. 
     \end{align}
     \end{subequations}

  %!!!!!!!!!!!!!!!!! - 6 - !!!!!!!!!!!!!!!!!!!!!!!!
  \item Coordinate translation  and parameter reflexion [from Eq.~\ref{sym_5}], 
     \begin{subequations}
     \label{shift1}
     \begin{align}
     {x}_n({x}_0+1/2,{y}_0+1/2; -\kappa_1,\kappa_2) &=
           {x}_n({x}_0,{y}_0;\kappa_1,\kappa_2) \nonumber \\ 
           & \,\,\,\,\, +1/2 + n \, , 
      \label{shift11}\\
      {y}_n({x}_0+1/2,{y}_0+1/2; -\kappa_1,\kappa_2) &= 
          {y}_n({x}_0,{y}_0;\kappa_1,\kappa_2) +1/2 \, .
      \label{shift12}
     \end{align}
     \end{subequations}
 \end{enumerate}
% =======================================
 Note that property $P3$ in (\ref{sym_3})  implies that if  there is an invariant circle above the line $y=0.5$, then  there is an invariant circle corresponding to its \emph{reflected} image bellow $y=0.5$.
This same property exists in the standard map, so its  invariant circle $\gamma$ has a reflected image (with rotation number $1-\gamma$) in the lower half of the cell, and both break up for the same value of the parameter $\kappa_G=0.971635406$.
 
%////////////////////////////
%////////////////////////////

%--------------------------------------------------
\subsection{Rotation number symmetries.}
\label{sec:RotationSym}
%--------------------------------------------------
 We define the rotation number of an orbit of the map (\ref{mapaut})
with parameters $(\kappa_1,\kappa_2)$, 
\bq
\omega ({x}_0,{y}_0; \kappa_1,\kappa_2)= \lim_{n \rightarrow \infty} \frac{{x}_n({x}_0,{y}_0; \kappa_1,\kappa_2)-{x}_0}{n} \, , 
\label{ome_r}
\eq
whenever the limit exist.

 Then from Eqs.~(\ref{shift31}), (\ref{shift21}) and (\ref{shift11}), respectively, it follows that 
\begin{eqnarray}
  \omega({x}_0+1/2,{y}_0; -\kappa_1,-\kappa_2)&=&
     {\omega}({x}_0,{y}_0;\kappa_1,\kappa_2)  
        \, , \nonumber \\ \label{ome1} \\
  \omega({x}_0,{y}_0+1/2; \kappa_1,-\kappa_2)&=&
     {\omega}({x}_0,{y}_0;\kappa_1,\kappa_2) 
    +1  \, , \nonumber \\ \label{ome2}\\
 \omega({x}_0+1/2,{y}_0+1/2; -\kappa_1,\kappa_2)&=&
     {\omega}({x}_0,{y}_0;\kappa_1,\kappa_2) 
    +1   \, . \nonumber \\ \label{ome3}. 
\end{eqnarray}

 Therefore, for each invariant circle passing through $({x}_0,{y}_0)$ with rotation number $\omega$ in the NASM with parameters $(\kappa_1,\kappa_2)$,  there exist  up to three other associated  invariant circles:   
one  with rotation number $\omega+1$  for $(-\kappa_1,\kappa_2)$ passing through $({x}_0+1/2,{y}_0+1/2)$,
one also with rotation number $\omega+1$ for $(\kappa_1,-\kappa_2)$ passing through $({x}_0,{y}_0+1/2)$,
and 
one with rotation number $\omega$ for $(-\kappa_1,-\kappa_2)$ passing through $({x}_0+1/2,{y}_0)$.

 In addition, for a given $(\kappa_1,\kappa_2)$, from Eq.~(\ref{int}) it can be shown that,
\begin{equation}
    \omega(1-{x}_0,1-{y}_0; \kappa_1,\kappa_2)=
      2 - {\omega}({x}_0,{y}_0;\kappa_1,\kappa_2)  \, , \label{ome4} \\
\end{equation}
 which is 
 a property that also applies to
the standard map (\ref{stdmap}), but with a different shift\footnote{$1$ instead of $2$ on the right hand side of Eq.~(\ref{ome4}).}.

Finally from Eq.~(\ref{sim_M}) and the definition of rotation number in (\ref{ome_r}), it follows that,
\begin{equation}
  \omega(\mathcal{S}_{\kappa_2}(\mathbf{z}_0);\kappa_2,\kappa_1)=
     \omega(\mathbf{z}_0;\kappa_1,\kappa_2) \,,
     \label{ome5}
\end{equation}
where $\mathbf{z}_0=(x_0,y_0)^T$.

%||||||||||||||||||||||||||||||||||||||||||||||||||
\section{Periodic orbits.}
%||||||||||||||||||||||||||||||||||||||||||||||||||
\label{sec:Periodic}
The periodic orbits are sets that in many cases offer information that can be used to characterize maps, approximate invariant sets and study linear stability. For these reasons,  it is important to study them in the case of the NASM.

The close relationship between $\mathcal{T}_{\kappa_1 \kappa_2}$ and $\mathcal{T}_{\kappa_2 \kappa_1}$ reflects on the fact that if $\mathbf{z}=(x,y)^T$ is an $n$-periodic orbit, on the lift of the map\footnote{To the universal cover of $\mathbb{T}\times\mathbb{R}$.} 
with rotation number $m/n$ of ${\cal T}_{\kappa_1 \kappa_2}$, that is,
\bq
 {\cal T}_{\kappa_1 \kappa_2}^n ({\bf z}) = {\bf z} + (m,0)^T \, ,
\eq
then from Eq.~(\ref{sim_M}) it follows that, 
\bq
{\bf w}= \mathcal{S}_{\kappa_2}^{-1} {\bf z} \, , \label{wz_property}
\eq
is an $m/n$-periodic orbit of ${\cal T}_{\kappa_2 \kappa_1}$, i.e.
\bq
  {\cal T}_{\kappa_2 \kappa_1} ^n ({\bf w}) = {\bf w} + (m,0)^T \, .
\eq

Note that the linear stability properties of $\mathbf{z}$ and $\mathbf{w}$  are the same because the trace of a product of matrices is invariant under the product commutation.

%|-|-|-|-|-|-|-|-|-|-|-|-|-|-|-|-|-|-|-|-|-|
\subsection{Period-one orbits}
\begin{figure}
\centering
\includegraphics[width=\linewidth, trim={1 1 0 0}]{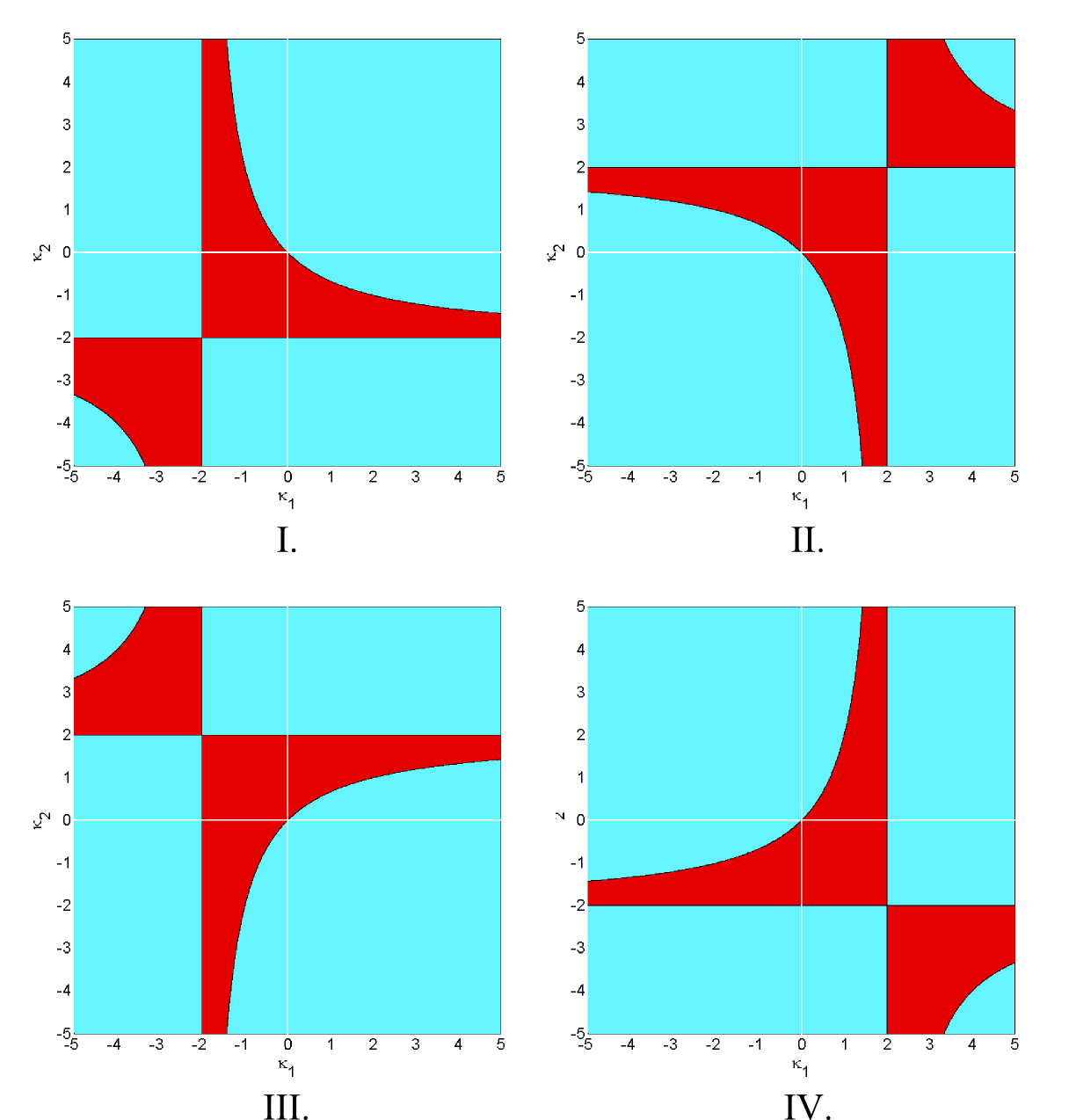}
\caption{Stability of primary period-one fixed points. Red regions correspond to stable orbits.
As expected, the stability diagrams satisfy the symmetry relations in Eqs.~(\ref{sym_1})-(\ref{sym_4}).}
\label{fig_stab_per_1}
\end{figure}

%|-|-|-|-|-|-|-|-|-|-|-|-|-|-|-|-|-|-|-|-|-|
There are six primary period-one orbits. By primary we mean that they exist for any values of $\kappa_1$ and $\kappa_2$. In addition there are bifurcated period-one orbits that exist only for certain values of $\kappa_1$ and $\kappa_2$. The primary orbits are
\bq
{\cal P}_1=\left\{ \left(1/2,0\right), \left(0,0\right), \left(0,\pm 1/2\right),  \left(1/2,\pm 1/2\right) \right\} \, .
\eq
The stability of these orbits is determined by the residue
\bq
R= \frac{1}{4} \left [ 2 -{\rm Tr} \left( \nabla \mathcal{T}_{\kappa_1 \kappa_2} \right) \right ] \, ,\label{greene_res}
\eq
where ${\rm Tr}$ denotes the trace, and $\nabla \mathcal{T}_{\kappa_1 \kappa_2}$ is the derivative of the map evaluated at the fixed point. A fixed point $(x_*,y_*)$ is stable if and only if $0<R<1$. From this it follows that,
\begin{eqnarray}
\label{period_one_stab}
\mathrm{I}.   &  
     \left(0,0\right)  & {\rm is \,\, stable\,\, iff\,\,}  
     0< - \kappa_1-\kappa_2 - \frac{\kappa_1 \kappa_2}{2} <2\,,
      \phantom{OO} \label{period_one_r1} \\ 
\mathrm{II}.  & 
     \left(\frac{1}{2},0\right)  & {\rm is \,\, stable\,\, iff\,\,}  
     0< \kappa_1+\kappa_2 - \frac{\kappa_1 \kappa_2}{2} <2\,, 
       \label{period_one_r2}\\
\mathrm{III}. & 
     \left(0,\pm \frac{1}{2} \right)  & {\rm is \,\, stable\,\, iff\,\,}  
     0<  \kappa_2-\kappa_1 + \frac{\kappa_1 \kappa_2}{2} <2\,,
       \label{period_one_r3}\\
\mathrm{IV}.   & 
     \left(\frac{1}{2},\pm \frac{1}{2} \right)  & {\rm is \,\, stable\,\, iff\,\,} 
     0<  \kappa_1-\kappa_2 + \frac{\kappa_1 \kappa_2}{2} <2\,. \label{period_one_r4}
\end{eqnarray}
Figure~\ref{fig_stab_per_1} shows the stability regions of the primary period-one fixed  points in the $(\kappa_1,\kappa_2)$ space, according to (\ref{period_one_r1})-(\ref{period_one_r4}). 
As expected, the results in 
Fig.~\ref{fig_stab_per_1} are consistent with the symmetries in(\ref{sym_1})-(\ref{sym_4})).

%%%%% FIGURE %%%%%%%%%%%%%%%%%%
%%%%%%%%%%%%%%%%%%%%%%%%%%%

The secondary period-one orbits $(x_*,y_*)$ are solutions to the system,
  \begin{subequations}
  \label{second}
  \begin{align}
   2y_* &= -\frac{\kappa_1}{2\pi}\sin (2\pi x_*) \qquad \mathrm{mod}(1)\,,  \label{second1}\\
   \kappa_1\sin (2\pi x_*) &= -\kappa_2\sin \left( 2\pi x_* +\frac{\kappa_1}{2}\sin (2\pi x_*)\right)\,. \label{second2}
  \end{align}
  \end{subequations}

\begin{figure}[h!]
     \centering 
     \includegraphics[width=\linewidth, trim={0.3cm 1.2cm 0 0}]{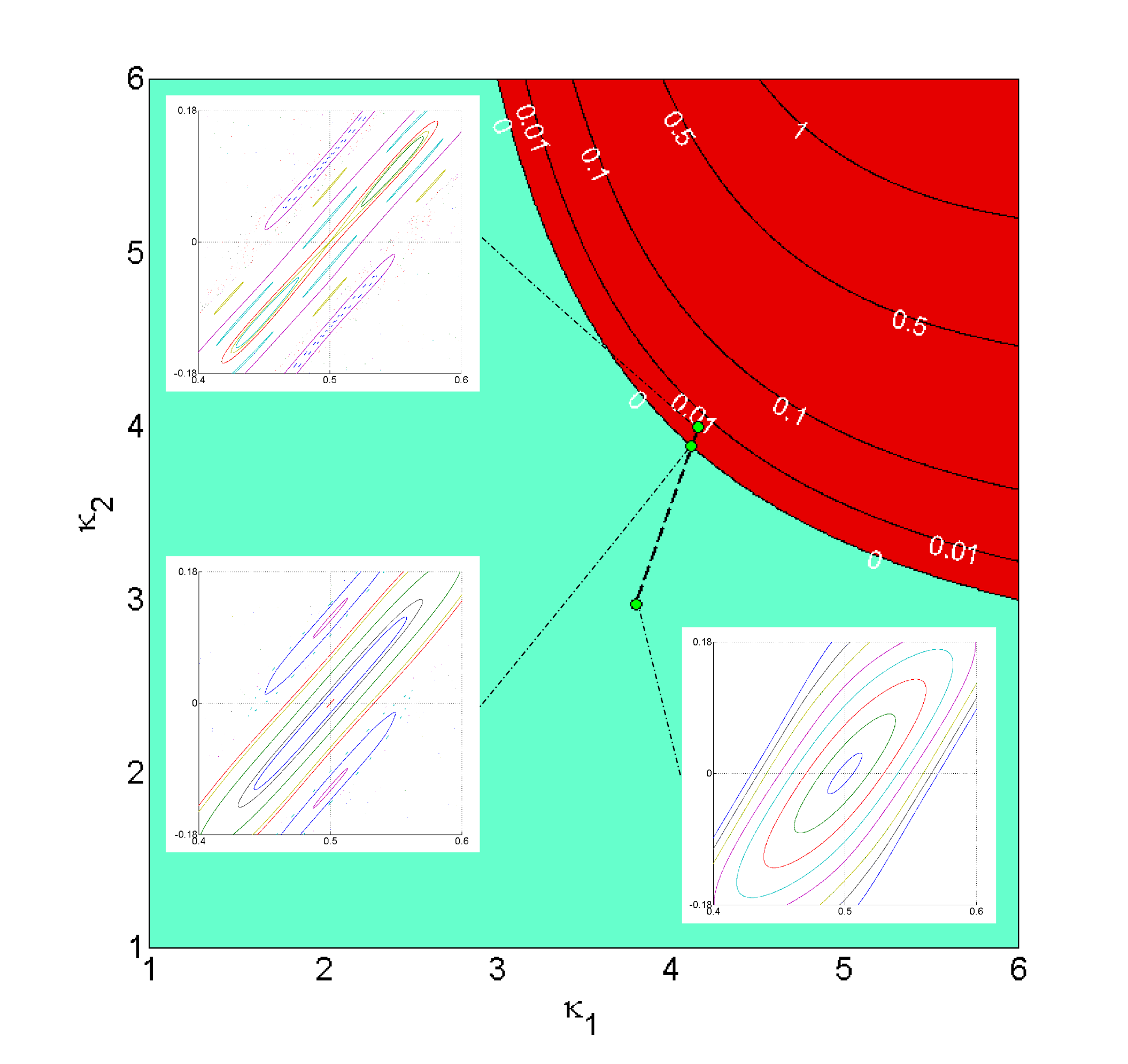} 
     \caption{Region of existence of the secondary period 1 periodic orbits around the elliptic point $\left(\frac{1}{2},0\right)$ in the positive quadrant of the parameter space. The light green shaded region corresponds to values of the parameters where the secondary periodic orbits do not exist. The figure displays in red the increasing values of the absolute error of the approximate value of $\mathsf{z}_*$ from (\ref{z_aprox}), when substituted into (\ref{second2}). The red colored region also coincides with values of the Greene's residue, Eq.~(\ref{greene_res}), for when the orbits are stable.    }
     \label{perio1s}
 \end{figure}

 A Taylor expansion of (\ref{second2}) around the elliptic point $(1/2,0)$, neglecting fifth order terms, around $|z_*|=|2\pi x_*-\pi|<<1$, allows to estimate $z_*$ as
\begin{equation}
 z_*^2= \frac{6\kappa_1+6\kappa_2+3\kappa_1\kappa_2}
       {\kappa_1-\frac{\kappa_1\kappa_2}{2}+\kappa_2\left( 1-\frac{\kappa_1}{2}\right)^3}
       \label{z_aprox}
\end{equation}
 Figure {\ref{perio1s}} shows the region in the positive quadrant of the parameter space where it is possible to find these associated periodic orbits of period 1 around the primary period-one fixed point $\left(\frac{1}{2},0\right)$. It should be noted that the limiting curve in the figure coincides with the limiting curve in Fig. \ref{fig_stab_per_1} for the same point.  These secondary families of periodic orbits appear after a pitchfork bifurcation of the primary fixed point $\left(\frac{1}{2},0\right)$.  The value of Greene's residuee evaluated at the orbits correspond to the region shaded in red tells us that these orbits are stable.

%////////////////////////////
%////////////////////////////

%||||||||||||||||||||||||||||||||||||||||||||||||||
\section{Transport barriers: known limit cases and direct computation.}\label{sec:KnownLimits}
%||||||||||||||||||||||||||||||||||||||||||||||||||
 Our goal is to find the critical parameter values  of the map (\ref{mapaut}) for which the map exhibits \emph{global transport} occurs.  We will say that there is global transport when there exists at least one initial condition $(x_0,y_0)$ such that the $y$ variable is unbounded.  It is established\cite{Greene79} that for this kind of map the only barriers to global transport are the invariant circles in the cylinder  that are not homotopic to a point. 
The map (\ref{mapaut}) is symplectic (Eq.~(\ref{symplectic})), exact (Eq.~(\ref{exact})) and twist (Eq.~(\ref{twist_cond})) for small values of the parameters. Because of these properties, Moser's Twist Theorem 
 guarantees  the existence  
invariant circles in a neighborhood of the origin in the parameter space, where 
 the map is integrable.
Then it is reasonable to 
compute the critical boundaries ($CB$) in the 
parameter space for the onset of global transport $ (CB_{gt})$ and for the destruction of invariant circles with a fixed rotation number, $\omega$, ($CB_{\omega}$). 
By $CB$ we mean the boundary of the open region ${\cal R}$ in the parameter space (assumed to be simply connected) such that for $(\kappa_1,\kappa_2)$ values in $\cal{R}$ there is not global transport in the case of ($CB_{gc}$)  or a given invariant circle with rotation number $\omega$ exists in the case of ($CB_\omega$).

We propose three different approaches to find these transport barriers: analytical reductions of the map for particular values of the parameters, direct numerical iteration for a range of values of the parameters  which gives upper bounds, and a continuation method over invariant circles with an a priori chosen rotation number in the parameter space  which gives lower bounds to global transport.

%|-|-|-|-|-|-|-|-|-|-|-|-|-|-|-|-|-|-|-|-|-|
\subsection{Known limit cases.}
%|-|-|-|-|-|-|-|-|-|-|-|-|-|-|-|-|-|-|-|-|-|
\label{subsec:Knownlim}
 In the standard map, the term \emph{critical parameter value}, $\kappa_c$,  
 is established in the literature\cite{Greene79} as the value for with any further increment $|\kappa|>|\kappa_c|$ there are no invariant curves. We call \emph{critical invariant curve} to such invariant curve that exist for $\kappa=\kappa_c$ and cease to exist after the critical value. It is a known conjecture\cite{Greene79} that for the standard map, the rotation number of the critical invariant curve is equal to the golden mean $\gamma=\frac{\sqrt{5}-1}{2}$. 

 We present the known reductions of map (\ref{mapaut}).
%########################################################
%########################################################
\begin{enumerate}
   %!!!!!!!!!!!!!!!!!!! - 1 - !!!!!!!!!!!!!!!!!!!!!!
    \item $(\kappa_1,\kappa_2)=(\kappa_1,0)$, $\kappa_1\geq0$.
    
     In this case the map (\ref{mapaut}) reduces to,
    \begin{subequations}
    \label{mapk1v1_old}
    \begin{align}
            x_{n+1}&= x_n + 2y_n + \frac{2\kappa_1}{2\pi}               
                       \sin(2\pi x_n) \label{mapk11v1}\\
            y_{n+1}&= y_n + \frac{\kappa_1}{2\pi} 
                                   \sin(2\pi x_n) \,  \label{mapk12v1} 
    \end{align}
    \end{subequations}               
which upon the change of coordinates: $\{X=x,Y=2y\}$,   becomes the standard map (\ref{stdmap}) with $\epsilon=2\kappa_1$.
 And it is well known in this case that the critical invariant circle 
has rotation number equal to the golden mean $\gamma$ and breaks at the critical value 
$\kappa_G\approx 0.971635406$ \cite{Greene79,Meiss92}. This implies that the {critical invariant circle} for the case
 $(\kappa_1,\kappa_2)=(\kappa_1,0)$ has rotation number $\omega_{c}=\gamma$ and breaks for $\kappa_1{}_c=\kappa_G/2$.  Applying  Eq.~(\ref{ome4}) to the standard map reduction and the NASM, it follows that there are also three more critical invariant circles  with rotation numbers: $1-\gamma$, 
 $2-\gamma$ and $\gamma+1$.

  %!!!!!!!!!!!!!!!!!!!!! - 2 - !!!!!!!!!!!!!!!!!!!!
  \item $(\kappa_1,\kappa_2)=(0,\kappa_2)$, $\kappa_2\geq0$.
  
     In this case the map (\ref{mapaut}) reduces to,
    \begin{subequations}
    \label{mapk2v1_old}
     \begin{align}
      x_{n+1} &= x_n + 2y_n + \frac{\kappa_2}{2\pi} 
                       \sin(2\pi\left( x_n +y_n\right)), \label{mapk2v11}\\
      y_{n+1}&= y_n  + \frac{\kappa_2}{2\pi} 
                       \sin(2\pi\left( x_n +y_n\right)).
                      \label{mapk2v12} 
     \end{align}
     \end{subequations}
 As before, with the change of coordinates: $\{X=x+y,Y=2y\}$, the map reduces to the standard map (\ref{stdmap})
  with perturbation parameter $\epsilon=2\kappa_2$.  Therefore the {critical invariant circle} has rotation number $\omega_c=\gamma$ and breaks for $\kappa_2{}_c=\kappa_G/2$.   From Eq.~(\ref{ome4}), there are three more critical invariant circles: $1-\gamma$ , $2-\gamma$ and $\gamma+1$. This case can also be considered a consequence of the Eq.~(\ref{ome5}) applied to the previous case.

    %!!!!!!!!!!!!!!!!!!!!! - 3 - !!!!!!!!!!!!!!!!!!!!
    \item $(\kappa_1,\kappa_2)=(\kappa,\kappa)$, $\kappa\geq0$.
    
      In this case the map reduces to 
the standard map iterated twice, which  that
the critical invariant circle of map (\ref{mapaut}), has rotation number $\omega_c=2\gamma$ and breaks at $\kappa_c= \kappa_G$. 
 By the symmetry (\ref{ome4}), there exists also the invariant circle  with rotation number: $\omega_{c_2}=2-2\gamma$.
     
Furthermore, applying the results of section {III.B} 
we can give the rotation numbers of the critical circles in the other quadrants. 
For example, for $(-\kappa,\kappa)$, the critical invariant circles are: $\omega_{c_1}=2\gamma-1$ and  $\omega_{c_2}=3-2\gamma$ for $\kappa_c=\epsilon_c$.

\end{enumerate}
%########################################################
%########################################################

It should  be noted that the change of variables used in the first two cases are  homotopic to the identity, so the barriers in the cylinder for the standard map reductions are barriers for the NASM as well.

%|-|-|-|-|-|-|-|-|-|-|-|-|-|-|-|-|-|-|-|-|-|
\subsection{Direct computation.}
%|-|-|-|-|-|-|-|-|-|-|-|-|-|-|-|-|-|-|-|-|-|
\label{subsec:direct}

We formally define a barrier to global transport as an invariant circle not homotopic to a point, which geometrically are circles that go around the cylinder $\mathbb{T}\times\mathbb{R}$. An invariant circle  of this kind, when it exists, always divides the phase space in two unbounded invariant regions, due to its invariance and the continuity of the map. 

 The numerical computation of the threshold to global transport is based on the following criteria:

\begin{proposition}
(Global transport criterion) \\
If for a set value of the parameters $(\kappa_1, \kappa_2)$, there is an initial condition $(x_0,y_0)$ with
 $y_0 \in   (0,1)$, for which  
$|y_n - y_0| >  2$ 
for some $n$ and all the invariant circles of the map have an amplitude\footnote{Where $\mathcal{A}$ is the difference between the highest and lowest point of the invariant circle.} $\mathcal{A}<1$, then the map has global transport.
\label{remak_gtc}
\end{proposition}
\emph{Proof.}
Lets assume that for the map with $(\kappa_1, \kappa_2)$, there exists an initial condition $(x_0,y_0)$ such that $|y_n-y_0|>2$ for a certain $n$ but there still exists an invariant circle fully contained in the cell $[0,1]\times(Q,Q+1)$ for $Q\in\mathbb{R}$. By the orbit symmetry property in Eq.~(\ref{trans}), using $r=0$ and $s=-\lfloor Q\rfloor:=-\mathrm{max}\{m\in\mathbb{Z}| m\leq Q\}$, there exists a \emph{copy} of the invariant circle in $[0,1]\times[0,2]$. Then, either the invariant circle or one of its \emph{copies} ($s=-\lfloor Q\rfloor\pm 1$) lies between $y_0$ and $y_n$. Which is a contradiction because $y_0$ and $y_n$  must be in the same connected component of the cylinder.

\hfill $\square$

Numerical evidence,  see e.g. Fig. \ref{cell_sym}(a), shows that at least for parameter values $|\kappa_i|<1$, the hypothesis, $\mathcal{A}<1$, holds.  The hypothesis is also obeyed because the invariant stable and unstable manifolds of the hyperbolic fixed point around $(0,0)$ 
obstruct the path of the invariant curve, and so 
do its integer translates in the $y$-direction.

\begin{figure}[h!]
  \centering
  \includegraphics[width=\linewidth,trim={0 1cm 0 0}]{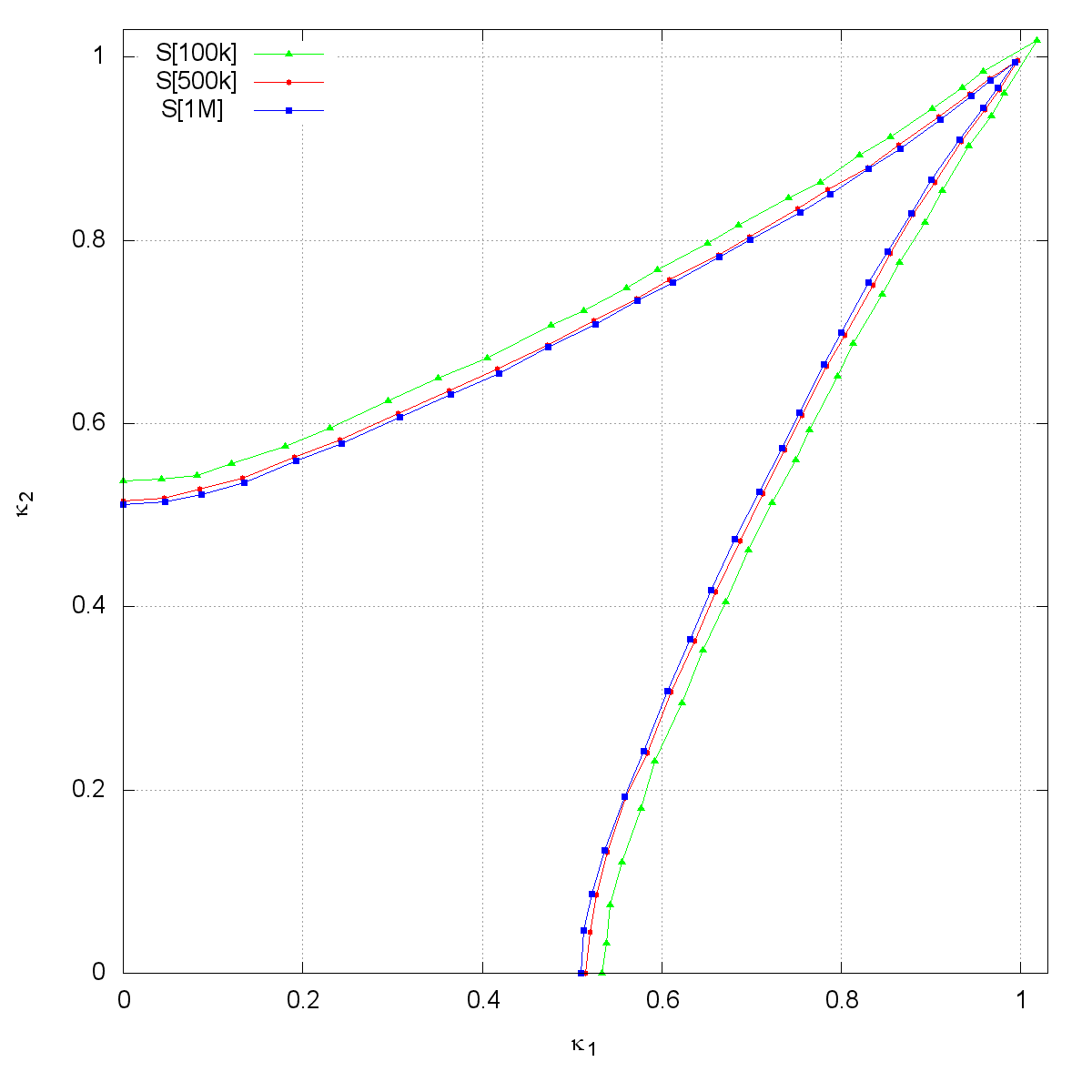} 
 \caption{ Critical boundary for global transport ($\mathit{CB}_{gt}$) in the non-autonomous  standard map in Eq.~(\ref{mapaut}), in the positive quadrant.
 The area  outside the ``horn" corresponds to parameter values for which there is \emph{global transport}, in the sense that at least one of $10^4$ initial condition uniformly distributed on the rectangle $[0,1 ]\times[0, 0.3]$ exhibited a displacement with:  
 ${y> 0.3 + 2}$ or ${y< - 2}$  after  
 $10^5$ (green triangles),  $5\times10^5$ (red circles)
 and $10^6$ (blue squares)  
 iterations of the map.}
     \label{simulc1}
 \end{figure}

\begin{figure}[h!]
     \centering
     \includegraphics[width=\linewidth,trim={0 1cm 0 0}]{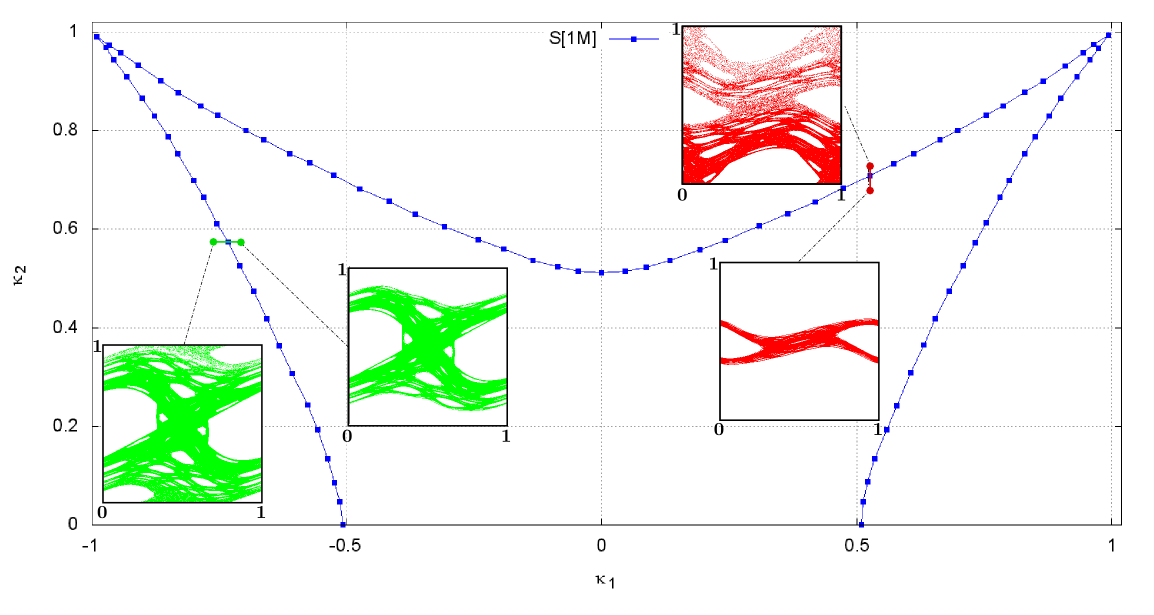} 
 \caption{ Critical boundary for global transport $(\mathit{CB}_{gt})$ in the non-autonomous  standard map in Eq.~(\ref{mapaut}) in the upper half plane of parameter space. As in Figure \ref{simulc1}, the area  outside the ``horns" corresponds to \emph{global transport}, for the same conditions as in Fig.~\ref{simulc1}.  The insets are images of the cell $[0,1]\times[0,1]$ of the phase space $(x,y)$ obtained from the iteration ($N=3\times10^6$) of the NASM for the initial condition $(0.5,0.44)$ for parameter values close the computed boundary: $(0.526,0.69)$, $(0.526,0.73)$, $(-0.69,0.525)$ and $(-0.73,0.525)$.}
     \label{simulc}
 \end{figure}

\begin{figure}[h!]
     \centering
     \includegraphics[width=\linewidth,trim={0 1cm 0 0}]{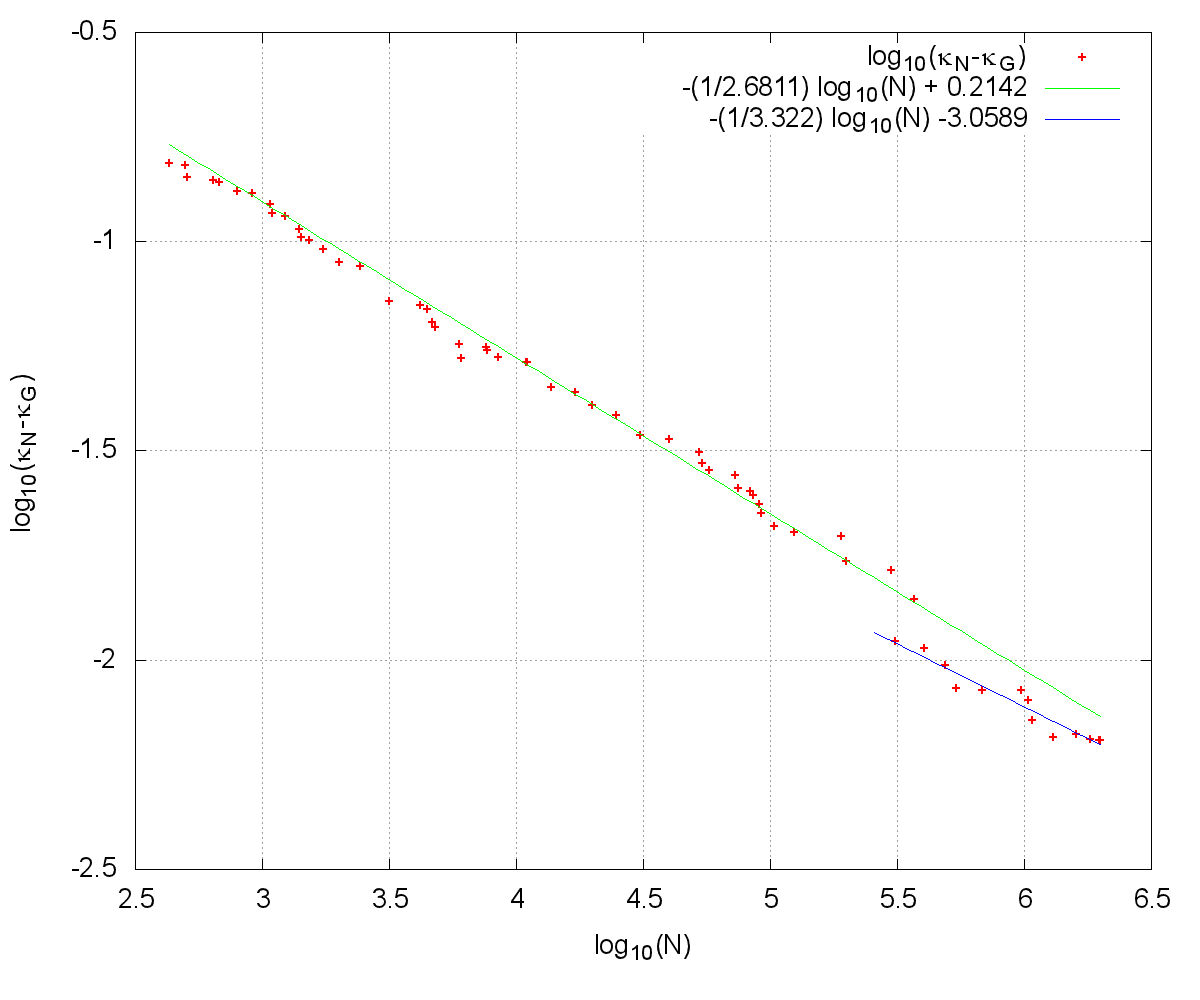} 
 \caption{ Convergence of the critical parameter value $\kappa_{c}$ as function of the number of iterates $N$, observed using the \emph{direct method} at the cusp of the \emph{horn}, along the diagonal $(\kappa_1,\kappa_2)=(\kappa,\kappa)$. To better display the results, the plotted quantity on the vertical axis is $\kappa_N-\kappa_G$.}
     \label{simul_conv}
 \end{figure}

We performed several series of $N$-iterations of the map (\ref{mapaut}) on $M$ initial conditions taken uniformly distributed in the rectangle $[0,1]\times[0,0.3]$ 
to determine, using the proposed transport criteria (Proposition 1) when the map displayed global transport for a wide range of the parameters $(\kappa_1,\kappa_2)$ inside the \emph{twist region}.  
Different number of initial conditions, $M$, were used,  finding $M=10^4$ to be a reasonable  tie between the computing capabilities and the consistency of the results. 
The procedure was repeated for an increasing number of iterations $N$. The convergence of the method is shown in Figure {\ref{simul_conv}}. 
 Figures \ref{simulc1} and \ref{simulc} show the  critical boundary for the global transport ($\mathit{CB}_{gt}$) found with these calculations, i.e. the locus of points in parameter space for which no critical invariant curves were detected.

In all the cases studied, the $\mathit{CB}_{gt}$ 
in the upper half plane of the parameter space	 were symmetric (up to machine precision) on the right (left) quadrant with respect to the line $\kappa_2=\kappa_1$ ($\kappa_2=-\kappa_1$), Fig. \ref{simulc}.
%////////////////////////////
%////////////////////////////
%////////////////////////////
%////////////////////////////

%////////////////////////////
%////////////////////////////

The slope of the adjusted line to the convergence data obtained with the \emph{direct method} for a variable number of iterations $N$, Fig.~\ref{simul_conv}, suggests that the convergence of the method is algebraic: $(\kappa_N-\kappa_G) \propto N^{-1/{\eta}}$, with $\eta=2.681$ (and $\eta=3.322$ for large $N$), which is close to the value,$\eta=3.012$, reported in literature\cite{Meiss92} for  transport in the standard map.

%||||||||||||||||||||||||||||||||||||||||||||||||||
\section{Transport barriers: continuation method.}
\label{sec:Continuation}
%|-|-|-|-|-|-|-|-|-|-|-|-|-|-|-|-|-|-|-|-|-|
In this section we describe the numerical implementation of the computation
and numerical continuation of invariant circles of the map in Eq.~(\ref{mapaut}).
We will omit most of the technical 
details and only discuss the tools that
we will use. For further mathematical details (e.g. function spaces, geometric preliminaries,
and Diophantine properties) the reader is referred to Ref.~\onlinecite{Lla-01, Cal-Lla-10,
Lla-Gon-Jor-Vil-05,Cal-Cel-Lla-13}.

\begin{definition}
\label{diophantine}
We say that the irrational number $\omega$ is Diophantine if for a given $\tau$ there is  a constant $\nu$ such that,
\[|\omega\cdot q - p| \geq \nu|q|^{-\tau},\,\,p \in \mathbb{Z},\,\,q\in\mathbb{Z}\setminus{\{0\}}.\]
\end{definition}
We will denote the set of all numbers satisfying Definition \ref{diophantine} by
$\mathcal{D}(\nu, \tau)$.

Then, we look for invariant circles of $\mathcal{T}_{\kappa_1 \kappa_2}$ on which the dynamics 
is conjugated to a rigid rotation by a fixed Diophantine rotation number $\omega$.

 The method we use is best understood in the constructive proof of the {KAM} theorem in Ref.~\onlinecite{Lla-Gon-Jor-Vil-05}, which
relies among other things
in a Newton iteration in the spirit of Nash-Moser theory,  see Ref.~\onlinecite{Zehnder75}.
For us, largest advantage of using the constructive proof 
is that it leads to a very efficient numerical algorithm.
We will make use of the Nash-Moser techniques to produce algorithms 
that will allow us to continue solutions
$K: \torus \to \mathbb{T}\times\mathbb{R}$ in 
a Banach space of smooth functions, to the following invariance 
equation  for a given twist map $F:\mathbb{T}\times\mathbb{R}\to \mathbb{T}\times\mathbb{R}$,
\begin{equation}\label{invariance}
F\circ K(\theta) = K(\theta+\omega),
\end{equation}
for $\omega \in \mathcal{D}(\nu, \tau)$.  Starting from the integrable case of the map $F$
and move the parameter as close to the breakdown of analyticity of the 
invariant circles as possible. We use the criterion of breakdown in Ref.~\onlinecite{Cal-Lla-10},
namely when we are close to the breakdown of
analyticity the derivatives of the solution $K$ start to blow up at points of
$K(\torus)$.

Continuation methods like the one presented here have already been used in
several contexts. See for instance, Ref.~\onlinecite{Calleja09, Cal-Lla-10b},
for models in statistical mechanics, Ref.~\onlinecite{Hug-Lla-Sir-06,
Fig-Luq-Har-15, book} for examples in symplectic maps, 
Ref.~\onlinecite{Cel-Cal-10, Cal-Fig-12} for conformally symplectic models, 
and Ref.~\onlinecite{Fox-Mei-13} 
for volume preserving maps.

The main idea of the method is to start form an approximate solution
of the invariance equation applied to $\mathcal{T}_{\kappa_1 \kappa_2}$. We will say that the solution $K_0$ is approximately
invariant if 
\[e_0(\theta) = \mathcal{T}_{\kappa_1 \kappa_2} \circ K_0 (\theta) - K_0(\theta + \omega),\]
and $\|e_0\|$ is a small function with respect to the norm $\|\cdot\|$ of the 
Banach space of smooth functions. 

We produce a ``better'' approximate solution (a solution that 
approximates the invariance equation \eqref{invariance} with
a smaller error), by adding a periodic function $\Delta: \torus\to \mathbb{R} \times \torus$ so that $K_1(\theta) = K_0(\theta) + \Delta(\theta)$
has an error
\[e_1(\theta) = \mathcal{T}_{\kappa_1 \kappa_2} \circ K_1 (\theta) - K_1(\theta + \omega),\]
with $\|e_1\| \approx \|e_0\|^2$.

Indeed, according to Nash-Moser theory, 
adding an appropriate correction $\Delta$ could provide an error satisfying
the quadratic property above. The correction $\Delta$ we could use, would solve the Newton step
equation,
\begin{equation}\label{newton}
 D\mathcal{T}_{\kappa_1 \kappa_2}(K_0(\theta))\Delta(\theta) - \Delta(\theta + \omega)
= - e_0(\theta).
\end{equation}
It is easy to check that if we were able to solve for $\Delta$ form equation
\eqref{newton}, then the norm of the new error, $\|e_1\|$ will be of order $\|e_0\|^2$.

One can try to solve numerically the newton equation in \eqref{newton}
for $\Delta$, but the most efficient methods
would require $O(n^2)$ operations, where 
$n$ is the number of points that one uses to
represent the invariant circle. The alternative that we will follow is to reduce
the Newton step equation \eqref{newton} by introducing a symplectic change of coordinates
around the approximate solutions. The implementation will yield methods that 
require $O(n\log n)$ operations. 

The change of coordinates is around an approximate solutions $K_0$ is given by
a $2$ by $2$ matrix composed of two column vectors. The first $2$ by $1$ column
is the vector $DK_0(\theta)$ representing the tangent bundle to the approximate
solution with base at every point $K(\theta)$. The second column vector is
a symplectic conjugate bundle which in our case of one dimensional circles
reduces to a vector orthogonal to the tangent bundle, namely 
$J^{-1}DK_0(\theta)N_0(\theta)$, where $J$ is the matrix representation of the symplectic form $\Omega$ and
\[N_0(\theta) {:}= [DK_0(\theta)^T DK_0(\theta)]^{-1}.\]
So the matrix can be written as follows,
\begin{equation}\label{Mmatrix}
M_0(\theta) = [DK_0(\theta) | J^{-1}DK_0(\theta)N_0(\theta)] \,.
\end{equation}

This change of coordinates is symplectic and transforms approximately the 
matrix $D\mathcal{T}_{\kappa_1 \kappa_2}(K_0(\theta))$ into an upper triangular 
matrix with ones along the diagonal, namely,
\begin{equation}\label{reduction}
D\mathcal{T}_{\kappa_1 \kappa_2}(K_0(\theta)) M_0(\theta) = M_0(\theta + \omega)
\left ( \begin{array}{cc}
1 & S_0(\theta) \\
0 & 1 \\
\end{array}
\right )
\end{equation}
where
\[S_0(\theta) = N_0(\theta + \omega) DK_0^T(\theta + \omega) 
D\mathcal{T}_{\kappa_1 \kappa_2}(K_0(\theta)) DK_0(\theta) N_0(\theta).\]
The function $S_0(\theta)$ is related to the local twist condition on the invariant
circle $K_0$.
We apply the change of coordinates to the Newton step, $\Delta(\theta) =
M_0(\theta) W(\theta)$, to reduce approximately the Newton step equation. 
Indeed, the new Newton step now has to solve  the equation,
\begin{eqnarray}\label{newtonreduced}
&&\left( \begin{array}{cc}
1 & S_0(\theta) \\
0 & 1 \\
\end{array} \right )
\left( \begin{array}{c}
W_1(\theta)\\
W_2(\theta)\\
\end{array} \right )  -
\left( \begin{array}{c}
W_1(\theta+\omega)\\
W_2(\theta+\omega)\\
\end{array} \right ) = \nonumber \\
&&  \phantom{OOOOOOOOOOOOOO} -M_0^{-1}(\theta + \omega) e_0(\theta).
\end{eqnarray}
Now, if we split the equation \eqref{newtonreduced} into components
we obtain two cohomological equations that we need to solve, namely
\begin{equation}\label{W2}
W_2(\theta) - W_2(\theta + \omega) = -[M_0^{-1}(\theta + \omega) e_0(\theta)]_2
\end{equation}
and
\begin{equation}\label{W1}
W_1(\theta) - W_1(\theta + \omega) = -[M_0^{-1}(\theta + \omega) e_0(\theta)]_1 - 
S(\theta) W_2(\theta).
\end{equation}
From equation \eqref{exact}, we know that
 $\mathcal{T}_{\kappa_1 \kappa_2}$ is exact. Then it can be shown  (\emph{Lemma 9} from Ref.~\onlinecite{Lla-Gon-Jor-Vil-05}) that 
\[\int_{\torus}[M_0^{-1}(\theta + \omega) e_0(\theta)]_2\,d\theta = {O}(\|e_0\|^2). \]
Therefore, we are able to solve for $\Delta$ form \eqref{newtonreduced} by
following the algorithm described below.
%@@@@@@@@@@@@@@@@@@@@@@@@@@@@@@@@@@@@@@@
\begin{algorithm}
  \label{algoritmo}
  \begin{itemize}
  \item[]
  \item[1)] Let $e_0 (\theta) = \mathcal{T}_{\kappa_1 \kappa_2} \circ K_0(\theta)
- K_0(\theta + \omega)$
  \item[3)] Compute the matrix $M(\theta)$ from equation \eqref{Mmatrix}.
  \item[4)] Solve for $W_2(\theta)$ from \eqref{W2}.
  \item[5)] Choose the average $\int_\torus W_2(\theta)d\theta$
   so that $-[M_0^{-1}(\theta + \omega) e_0(\theta)]_1 - S_0(\theta) W_2(\theta)$ has an average close to zero.
  \item[6)] Solve for $W_1(\theta)$ from \eqref{W1}.
  \item[7)] Compute the step $\Delta$,
\[\Delta(\theta) = M_0(\theta) W(\theta)\]
  \item[6)] Obtain the new parameterization $K_1$,
\[K_1(\theta) = K_0(\theta) + \Delta(\theta)\]
  \item[7)] Set $K_0(\theta) = K_1(\theta)$ and go to step 1).
  \end{itemize}
\end{algorithm} 
%@@@@@@@@@@@@@@@@@@@@@@@@@@@@@@@@@@@@@@@

\begin{remark}
One can verify that all the operations required to implement algorithm 
\ref{algoritmo}, are either diagonal in Fourier space of in real space.
To transform from real space to Fourier space one can  
use a Fast Fourier Transform (FFT), which is the most
expensive operation in the Algorithm \ref{algoritmo} in terms of 
arithmetic operations. Therefore, the cost of implementing
algorithm \ref{algoritmo} is $O(n\log n)$ operations,  where $n$ is the number of points used to represent the circle.
\end{remark}

In the following section we emphasize that there are useful relations
between the parameterization method exposed here and the symmetires 
that were introduced in Section \ref{sec:Symmetries}.

%---------------------------------------------------------
\subsection{Symmetries of the parameterization of an invariant circle}
%---------------------------------------------------------

In this Section we present some of the symmetries discussed
in Section \ref{sec:Symmetries} from the point of view of the parameterization
method. The goal is to  rewrite the symmetries
in terms of compositions of functions. In this way, we can use this composition 
formulation to rewrite the invariance equation that a parameterization 
of a circle with a certain symmetry should satisfy.  This way
we know that if a map has an invariant cicle, then we can use that parameterization
function, $K$, and transform it using the composition operators with respect to a given 
symmetry.  In this way, we construct a new parameterization that satisfies the invariance
equation of a new map with the given symmetry. We start by verifying that the invariant circles of 
the map ${\cal T}_{\kappa_1, \kappa_2}$ are expected to exist and to
have the rotation numbers found in Section \ref{subsec:Knownlim}. 
 
First, we consider that the map $F$ has an invariant circle
with rotation number $\omega$, whenever Eq.~(\ref{invariance}) is satisfied.
For instance, if $S_{\kappa}$  in Eq.~(\ref{stdmap}), has an invariant circle of rotation number $\omega$, with graph $K$,
then ${\cal T}_{\kappa \kappa}$ will have an invariant torus with rotation 
number $2 \omega$ as discribed above. Since, 
\[ S_{\kappa} \circ K(\theta) = K(\theta + \omega)\]
and
\[ S_{\kappa} \circ K(\theta + \omega) = K(\theta + 2\omega),\]
then
\begin{equation}\label{mapsquare}
 S_{\kappa} \circ S_{\kappa} \circ K(\theta) = S_{\kappa} \circ K(\theta + \omega)
= K(\theta + 2\omega).
\end{equation}
We also notice that if we define the following function for $\phi\in \mathbb{R}$,
\[R_\phi(x, y) = (x +\phi  \mod 1, y),\]
then the standard map \eqref{mapaut} satisfies that,
\begin{equation}\label{stdsymmetry}
S_\epsilon \circ R_{-{1\over 2}} = R_{-{1\over 2}} \circ S_{-\epsilon},
\end{equation} 
for any $\epsilon$.
We have one more symmetry of the parameterization that will be important.  
To see this, we write the parameterization of an invariant circle of 
a symplectic map of the cylinder, we write the components of the 
parameterization $K$ as follows. Let $u(\theta)$ be a $1$-periodic function of 
$\torus$, then 
\begin{equation}\label{param-periodic}
K(\theta) = \left(\begin{array}{c}
\theta + u(\theta) \\
\omega + u(\theta) - u(\theta - \omega)
\end{array}\right).
\end{equation}
It is clear from \eqref{param-periodic} and the 
periodicity of $u$, that 
\begin{equation}\label{circ-symmetry}
R_{-1} \circ K(\theta) = K(\theta -1).
\end{equation}

In particular, if $K$ is the parameterization of an invariant circle of 
$S_{\kappa}\circ R_{-{1\over 2}}$ then by the property 
\eqref{mapsquare} of the square of a map above, we have that
\[S_{\kappa}\circ R_{-{1 \over 2}} \circ S_{\kappa}\circ R_{-{1 \over 2}} 
\circ K(\theta) = K(\theta + 2 \omega).\]
So by the symmetry property in \eqref{stdsymmetry} of the Standard map 
together with the symmetry \eqref{circ-symmetry} of the invariant circle,
we have that 
\begin{eqnarray*}
&&S_{\kappa}\circ S_{-\kappa} \circ {R}_{-1} \circ K(\theta) = K(\theta + 2\omega) \\ 
&&\phantom{OOOOOOOO} \implies S_{\kappa}\circ S_{-\kappa} \circ K(\theta -1) = K(\theta + 2\omega)
\end{eqnarray*}
This is equivalent to saying that $K$ is an invariant circle
of $S_{\kappa}\circ S_{-\kappa}$ with rotation number $2\omega +1$, 
where $K$ is the invariant circle of $S_{\kappa}\circ R_{-{1\over 2}}$
with rotation number $\omega$,  which is the same as the symmetry in Eq.~(\ref{ome2}).

We verify the properties that were described in Section \ref{sec:KnownLimits}
were we anticipate that the maps $\mathcal{T}_{\kappa_1 0}$ and $\mathcal{T}_{0 \kappa_2}$
{can be rescaled} to standard maps with twice the perturbation parameter.

If we define the transformation
\[P_\xi(x, y) = (x \mod 1, \xi y),\]
as simple computation tells us that 
\[\mathcal{T}_{\kappa_1 0} = S_0 \circ S_{\kappa_1} = P_{1/2} \circ S_{2 \kappa_1 } \circ P_{2} \]
So if $S_{2\kappa_1}$ has an invariant circle, then the invariance equation is 
\[S_{2\kappa_1}\circ K(\theta) = K(\theta+\omega).\]
We immediately know that $P_{1/2} \circ K(\theta)$ is an invariant circle for 
$\mathcal{T}_{\kappa_1 0}$. Namely,
\[P_{2} \circ S_0 \circ S_{\kappa_1} \circ P_{1/2} \circ K(\theta) = K(\theta+\omega),\]
is equivalent to 
\[\mathcal{T}_{\kappa_1 0} \circ P_{1/2} \circ K(\theta) = P_{1/2} \circ K(\theta+\omega).\]
The case of the map $\mathcal{T}_{0 \kappa_1}$ is obtained for the above since when we use
\eqref{sim_M}, we have that $\mathcal{T}_{0 \kappa_1}= S_0^{-1}\circ\mathcal{T}_{\kappa_1 0}\circ S_0$.
Then, it is clear that the circle $S_0^{-1} \circ P_{1/2} \circ K$ is invariant for
$\mathcal{T}_{0 \kappa_1}$.  In fact, from Eq.~{\ref{sim_M}} it is easy to see that if $\mathcal{T}_{\kappa_1 \kappa_2}$ has an invariant circle with rotation number $\omega$, then $\mathcal{T}_{\kappa_2 \kappa_1}$ has an invariant circle with the same rotation number. This result was also stated in Eq.~(\ref{ome5}).

%-|-|-|-|-|-|-|-|-|-|-|-|-|-|-|-|-|-|-|-|
\subsection{Results from the parameterization method}
\begin{figure}[h!]
     \centering
     \includegraphics[width=\linewidth,trim={0 1cm 0 0}]{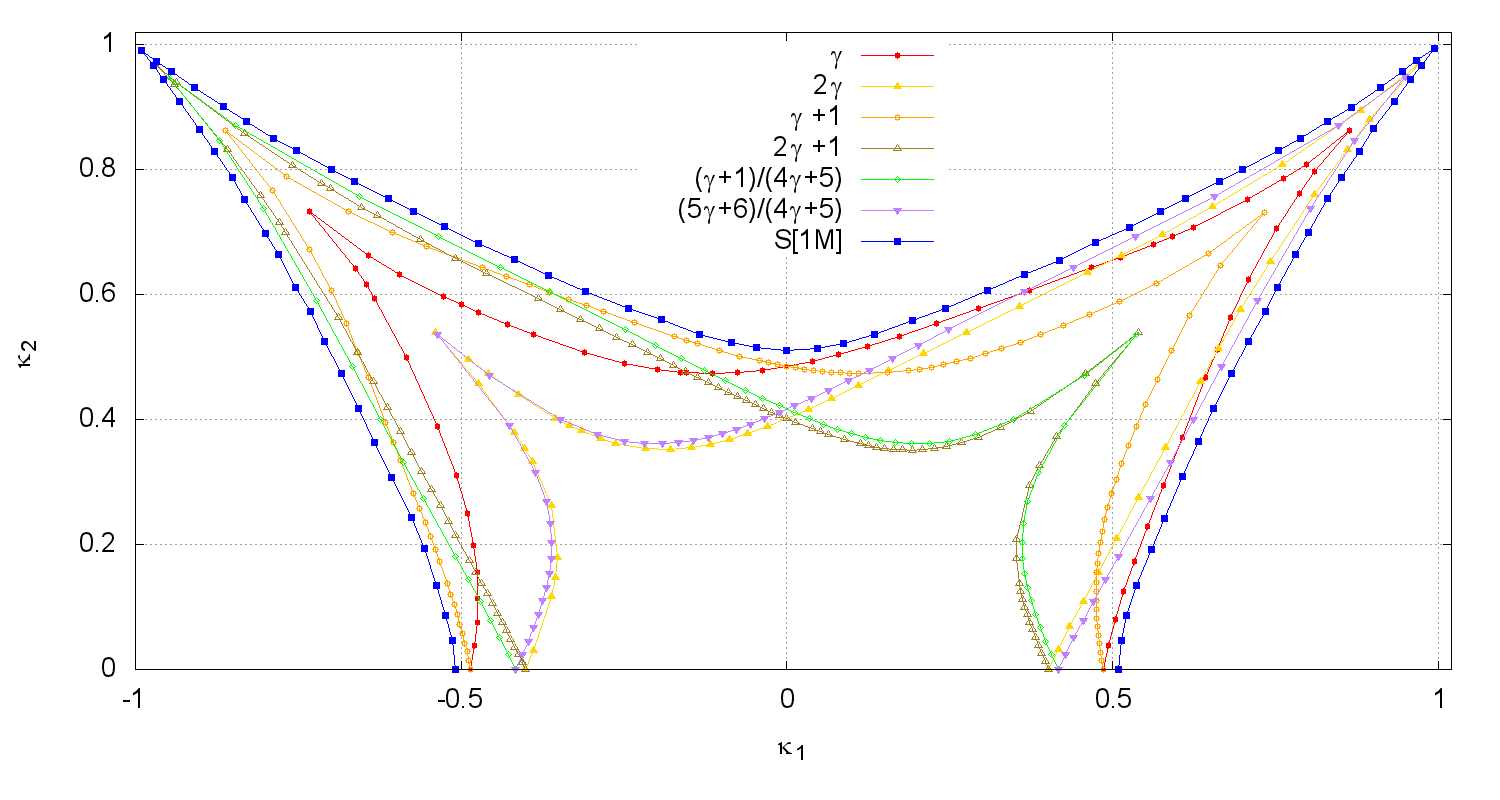} 
     \caption{ Critical boundaries for the existence of four KAM tori ($\mathit{CB}_\omega$) \{$\gamma$, $\gamma+1$, $2\gamma$, $2\gamma+1$, $\frac{5\gamma+6}{4\gamma +5}$ and $\frac{\gamma+1}{4\gamma +5}$\} found with the \emph{parameterization method} contrasted with the critical boundary ($\mathit{CB}_{gt}$) found via direct simulation of the NASM.  }
     \label{bif_fig}
 \end{figure}

%-|-|-|-|-|-|-|-|-|-|-|-|-|-|-|-|-|-|-|-|
\label{subsec:results_par}
 The parameterization method was applied using the rotation numbers that seemed relevant from section \ref{sec:KnownLimits}: $\gamma$, $\gamma+1$, $2\gamma$ and $2\gamma+1$ and also with $\frac{5\gamma+6}{4\gamma +5}$ and $\frac{\gamma+1}{4\gamma +5}$ which where found heuristically. 
The method 
yielded in all the cases as critical boundaries 
between the existence and \emph{breaking} of a given invariant circles with rotation number $\omega$, $(\mathit{CB}_\omega)$, 
a \emph{two horn}-shaped asymmetrical curve in the upper plane of the parameter space $(\kappa_1,\kappa_2)$.  The points obtained in the $\mathit{CB}_\omega$ correspond to the divergence of the norm of the derivatives of the parameterization function $K$. The critical boundaries are displayed in Figure \ref{bif_fig} with the addition of the a curve of the first method to contrast the results. All the  critical boundaries agree with the results from section \ref{sec:KnownLimits} and are contained in  the \emph{tightest} curve from  previous subsection. 
Notice that the $\mathit{CB}$ corresponding to $\gamma$ is related to the one of $\gamma+1$ and the $\mathit{CB}$ of $2\gamma$ is related to the one of $2\gamma+1$ by a reflection with respect to the $\kappa_2$-axis.

%////////////////////////////
%////////////////////////////

%////////////////////////////
 \begin{figure}[h!]
     \centering 
     \includegraphics[width=\linewidth,trim={0 1cm 0 0}]{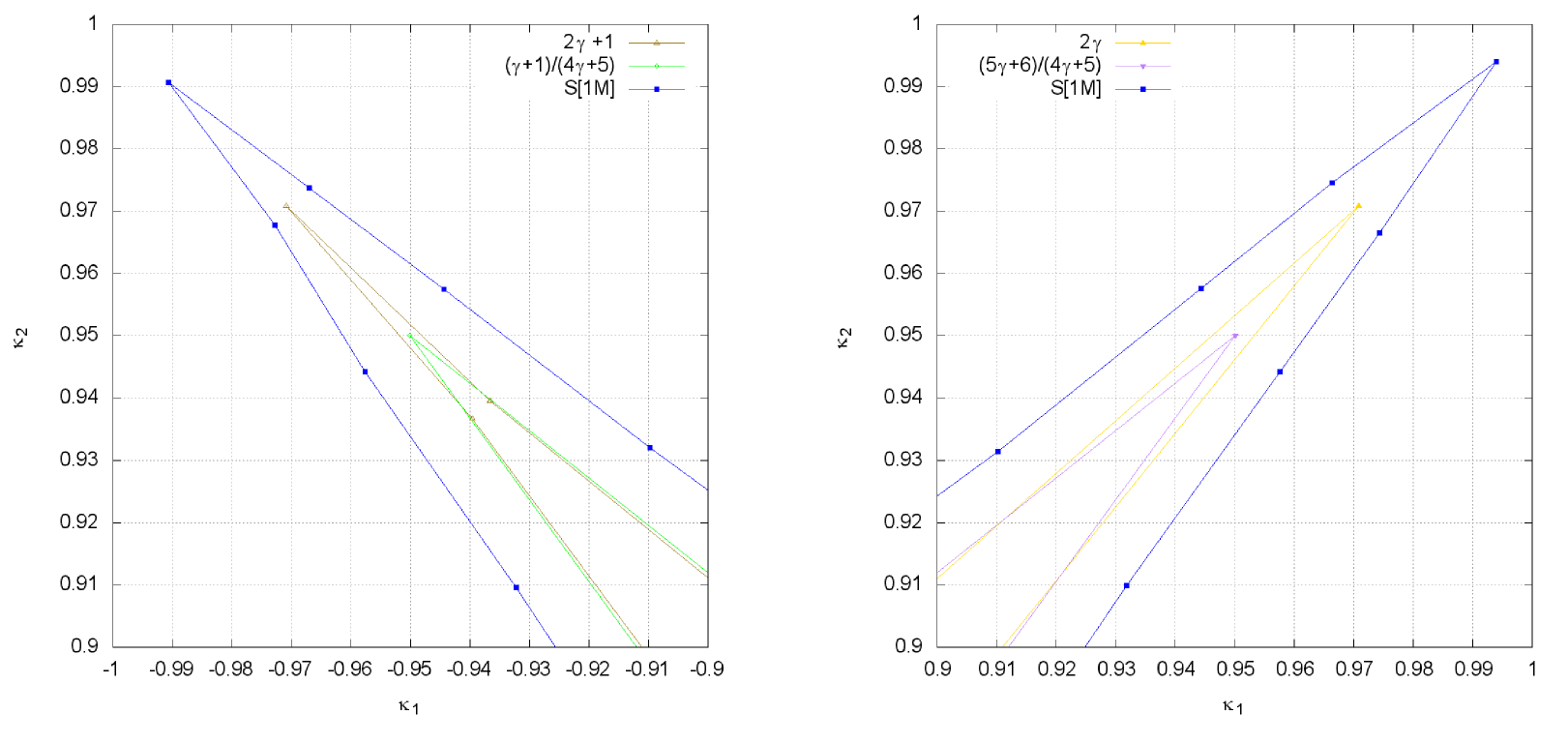}
 \caption{Details on the \emph{tips} of the \emph{left} and \emph{right} horns of Figure \ref{bif_fig}. }
     \label{bif_fig2}
 \end{figure}
%////////////////////////////

Additionally to the rotation numbers that were predicted to be of importance in section \ref{subsec:Knownlim}, the parameterization method was applied to two additional ones: $\frac{5\gamma+6}{4\gamma +5}$ and $\frac{\gamma+1}{4\gamma +5}$. 
 These number were obtained in an empirical form by iterating the map close to the invariant circle, obtaining a few digits by approximating the limit in (\ref{ome_r}) and then adding a tail of ones to the continued fraction.
Neither of these two numbers correspond to a known reduced case of the map (\ref{mapaut}), however the results displayed on Figure \ref{bif_fig} show that the corresponding $\mathit{CB}_\omega$'s give reasonable lower bounds to the $\mathit{CB}_{gt}$ found by direct method in some regions of the parameter space.

%||||||||||||||||||||||||||||||||||||||||||||||||||
\section{Discussion and conclusions}
%||||||||||||||||||||||||||||||||||||||||||||||||||
\label{sec:Discussion}
In section \ref{sec:KnownLimits} we found from the direct computation a horn shaped  
 critical boundary for global transport, $\mathit{CB}_{gt}$, giving the
threshold between bounded evolution of a set of initial conditions and global transport. Using the parameterization method in Sec.~\ref{sec:Continuation} we  found the critical boundaries for the existence of invariant curves with given rotation numbers, the $\mathit{CB}_{\omega}$'s. 
All the $\mathit{CB}_\omega$ were found to be fully contained inside the $\mathit{CB}_{gt}$. 
It is expected that the $\mathit{CB}_{gt}$ is  the convex hull of all the $\mathit{CB}_\omega$ associated to the invariant circles that exist for the map $\mathcal{T}_{\kappa_1 \kappa_2}$.

\begin{figure}[h!]
     \centering
     \includegraphics[width=\linewidth]{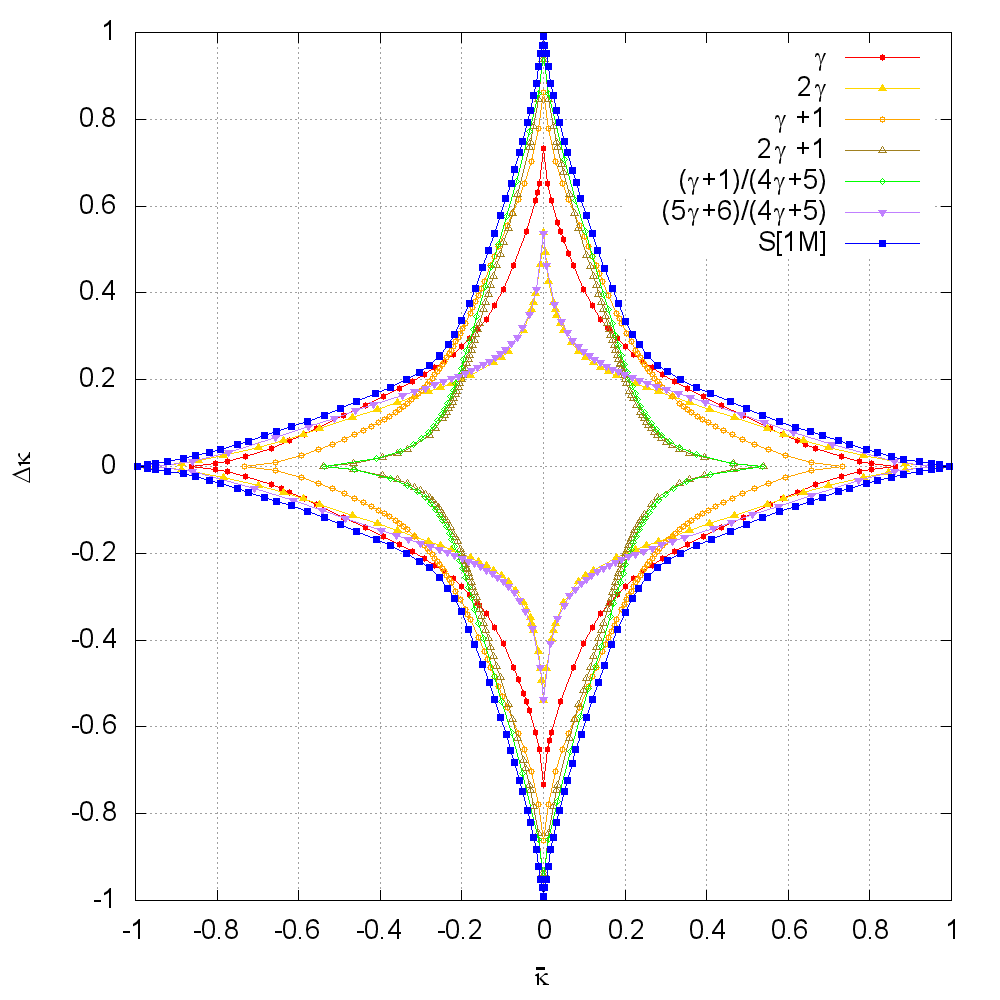} 
     \caption{ Plot in the $(\bar{\kappa},\Delta\kappa)$ parameter variables of the $\mathit{CB}_{\omega}$ of four KAM tori \{$\gamma$, $\gamma+1$, $2\gamma$, $2\gamma+1$, $\frac{5\gamma+6}{4\gamma +5}$ and $\frac{\gamma+1}{4\gamma +5}$\} found with the \emph{parameterization method} and the $\mathit{CB}_{gt}$ found via direct simulation of the NASM. The extension to the whole parameter space is obtained using the Eq.~(\ref{ome1}).   }
     \label{diamond1}
 \end{figure}

\begin{figure}[h!]
     \centering
     \includegraphics[width=12cm]{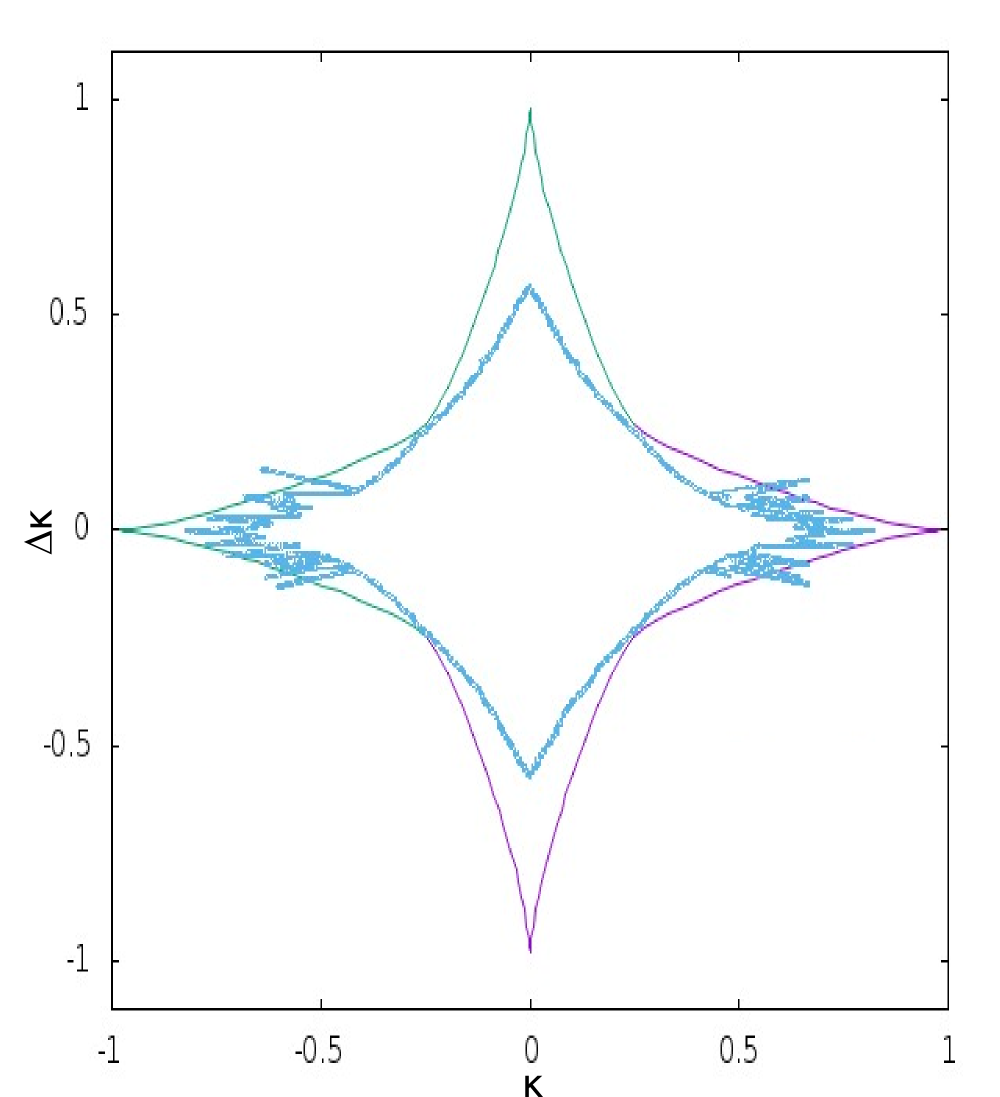} 
     \caption{ Critical boundary curve ($\mathit{CB}_{gt}$) found via direct simulation of the NASM compared with one of the $\mathit{CB}$ for a periodic orbit with rotation vector $(\frac{1705}{3136},\frac{2631}{3136})$ for map (\ref{rotmap})  taken directly from Ref.~\onlinecite{tompaidis96} [\emph{Exp. Math.} {\bf 5}(3):211-230,1996].}
     \label{diamond2}
 \end{figure}

$\mathit{CB}$'s were also computed using a Greene's residue method  
in Ref.~\onlinecite{tompaidis96} for the $3D$ map,
\begin{equation}
 \mathcal{F}_{k,\Delta k} 
 \left( \begin{array}{c}  x\\ y \\ \phi \end{array} \right)  =
    \left(\begin{array}{cl}
   {x}+ y+ \frac{k'}{2\pi} \sin(2 \pi x) & \mod 1\\
   y+ \frac{k'}{2\pi} \sin(2 \pi x) &\\ 
   \phi + \Omega  & \mod 1
\end{array} \right)\,,
\label{rotmap}
\end{equation}
 where $k'=\bar{\kappa}+\Delta\kappa\cos(2\pi\phi)$, which corresponds to a more general (quasi-periodic) variation of $\kappa$ than the one in Eq.~(\ref{kap}). 
 In the case $\Omega=1/2$, the map (\ref{rotmap}) is equivalent to the map in Eq.~(\ref{nas_map}) with the appropriate choice of parameters. Note that in the special case when the initial value of $\phi$ is $1/4$, both maps reduce to the standard map. When the initial value of $\phi$ is $0$,  the parameters of the maps are related by:$\bar{\kappa}=\frac{\kappa_1+\kappa_2}{2}$, $\Delta\kappa=\frac{\kappa_2-\kappa_1}{2}$.
In Ref.~\onlinecite{tompaidis96} the $\mathit{CB}$ displayed correspond to the critical values of the Greene's residue\cite{Greene79} in the parameter space $(\bar{\kappa},\Delta\kappa)$ for map (\ref{rotmap}) for different periodic orbits with rotation vectors that approximate a two dimensional invariant torus.
 These $\mathit{CB}$'s have a diamond shape similar to the curves obtained in Sec.~{\ref{sec:KnownLimits}} in $(\bar{\kappa},\Delta\kappa)$ parameter variables (see Fig.~\ref{diamond1}), which represent the threshold of global transport and disappearance of all invariant one dimensional tori for map (\ref{rotmap}) with rotation vectors of the form $(\omega,1/2)$. In figure \ref{diamond2} we compare our $\mathit{CB}_{gt}$ with one of the $\mathit{CB}$ from Ref.~\onlinecite{tompaidis96} corresponding to the periodic orbit with rotation vector $(\frac{1705}{3136},\frac{2631}{3136})$ with period $3136$, that approximates a golden rotation vector.
The critical curve is contained inside ours except for a few points that come out of the region we computed. 

To further explore this idea, we have studied map (\ref{rotmap}) with $\Omega=1/3$, which correspond to a $\kappa_n$ variation for map (\ref{nas_map}) with three values. Finding also critical boundaries that surpassed in a small region the one found for map (\ref{mapaut}).

%////////////////////////////
%////////////////////////////
%////////////////////////////
%////////////////////////////
Reference \onlinecite{simo16}  
computed the CB for a given rotation number in a \emph{driven standard map} similar to Eq.(53). As in the previous case, the results are consistent with ours. In particular, they also found a diamond shaped $\mathit{CB}$ in the ($\bar{\kappa}, \Delta \kappa$) parameter space. 

It is important to remark that in map (\ref{rotmap}) the two dimensional and the one dimensional tori are topological barriers to global transport due the uncoupled variation of $\phi$ with respect to $(x,y)$. In general, Fig.~\ref{diamond2} suggests that the destruction of the one dimensional tori implies the breaking of the two dimensional ones except for a few \emph{peaks} that stand outside of our curve in the parameter space. It remains to study why are the two dimensional tori more robust that the one dimensional ones in these \emph{peak} regions. Numerical evidence leads us to think that the critical boundary for global transport in the map (\ref{rotmap}) with $\Omega=p/q$ might fully contain the ones in Ref.~\onlinecite{tompaidis96}.

%|||||||||||||||||||||||||||||||||||||||||||||||||||||
\section*{Acknowledgments}
%|||||||||||||||||||||||||||||||||||||||||||||||||||||
This work was supported by PAPIIT IN104514, FENOMEC-UNAM and
by the Office of Fusion Energy Sciences of the US Department
of Energy at Oak Ridge National Laboratory, managed by UT-Battelle, LLC, for the
U.S.Department of Energy under contract DE-AC05-00OR22725.  We also express our gratitude  to LUCAR-IIMAS for making avaliable the CPU cluster and to the graduate program in Mathematics of UNAM for making the GPU servers available to perform our computations and especially to Ana Perez for her invaluable help. 

\bibliographystyle{unsrt}
\bibliography{NASM}

\end{document}